\documentclass[12pt]{article}
\usepackage{epsfig}
\usepackage{amsfonts}
\def\R{{\mathbb R}}
\def\C{{\mathbb C}}
\def\L{{\cal L}}
\def\Arf{{\it Arf}}
\def\sbp{\subparagraph}
\title{\hspace{\fill}{\normalsize\it To V.I.Arnold with admiration\vspace{15pt}}\\
\bf An integral generalization of the Gusein-Zade--Natanzon theorem}
\author{Sergei Chmutov}
\date{}
\begin{document}
 \maketitle

One of the useful methods of the singularity theory, the method of
real morsifications \cite{AC-0, GZ} (see also \cite{AGV}) reduces
the study of discrete topological invariants of a critical point
of a holomorphic function in two variables to the study of some
real plane curves immersed into a disk with only simple double
points of self-intersection. For the closed real immersed plane
curves V.I.Arnold \cite{Ar} found three simplest first order
invariants $J^\pm$ and $St$. Arnold's theory can be easily
adapted to the curves immersed into a disk. In \cite{GZN}
S.M.Gusein-Zade and S.M.Natanzon proved that the $\Arf$ invariant
of a singularity is equal to $J^-/2 (mod\ 2)$ of the corresponding
immersed curve. They used the definition of the $\Arf$ invariant
in terms of the Milnor lattice and the intersection form of the
singularity. However it is equal to the $\Arf$ invariant of the
link of singularity, the intersection of the singular complex
curve with a small sphere in $\C^2$ centered at the singular
point. In knot theory it is well known (see, for example
\cite{Ka}) that $\Arf$ invariant is the $mod\ 2$ reduction of some
integer-valued invariant, the second coefficient of the Conway
polynomial, or the Casson invariant \cite{PV}. Arnold's $J^-/2$ is
also an integer-valued invariant. So one might expect a relation
between these integral invariants.

A few years ago N.A'Campo \cite{AC-1, AC-2} invented a construction
of a link from a real curve immersed into a disk. In the case of
the curve originating from the real morsification method the link is
isotopic to the link of the corresponding singularity. But there
are some curves which do not occur in the singularity theory. 
In this article we
describe the Casson invariant of A'Campo's knots as a $J^\pm$-type 
invariant of the immersed curves. Thus we get an integral generalization 
of the Gusein-Zade--Natanzon theorem. 
It turns out that this $J^\pm_2$ invariant is a second order invariant 
of the mixed $J^+$- and $J^-$-types.
To the best of my knowledge, so
far nobody tried to study the mixed $J^\pm$-type invariants. It seems
that our invariant is one of the simplest such invariants. The 
problem of describing all second order $J^\pm$-type invariants is 
open.

In section {\bf 1} we review the A'Campo construction and list the
properties of the links obtained. In sections {\bf 2} and {\bf 3} we
introduce the Casson invariant and Arnold's type invariants of curves
immersed into a disk. In section {\bf 4} we formulate our main result.
The proof is based on M.Hirasawa's construction \cite{Hir} of a 
Seifert surface for the A'Campo links. We review this construction 
in section {\bf 5}.

\section{A'Campo divides and their links.}

\sbp{1.1. Definition} (\cite{AC-1, AC-2}). {\it A divide} $D$
is
the image of a generic immersion of a finite number of copies of the
unit interval $I$=[0,1] in the unit disk $B\subset\R^2$ such that
$\partial I\subset \partial B$. Here the word ``generic" means that
double points are the only singularities allowed and $D$ is
transversal to $\partial B$.

Similar object were known as {\it long curves} \cite{Ta, GZN}.

We consider divides up to isotopy of the disk $B$. The isotopy is
not assumed to be identical on the boundary $\partial B$.

\sbp{1.2. Example.} The curve $x^3+y^4=0$ has a singularity of type
$E_6$ at the origin \cite{AGV}. A small perturbation of it is a divide
which looks as follows.
$$
  \mbox{\begin{picture}(60,60)(0,0)
  \put(0,0){\epsfxsize=60pt \epsfbox{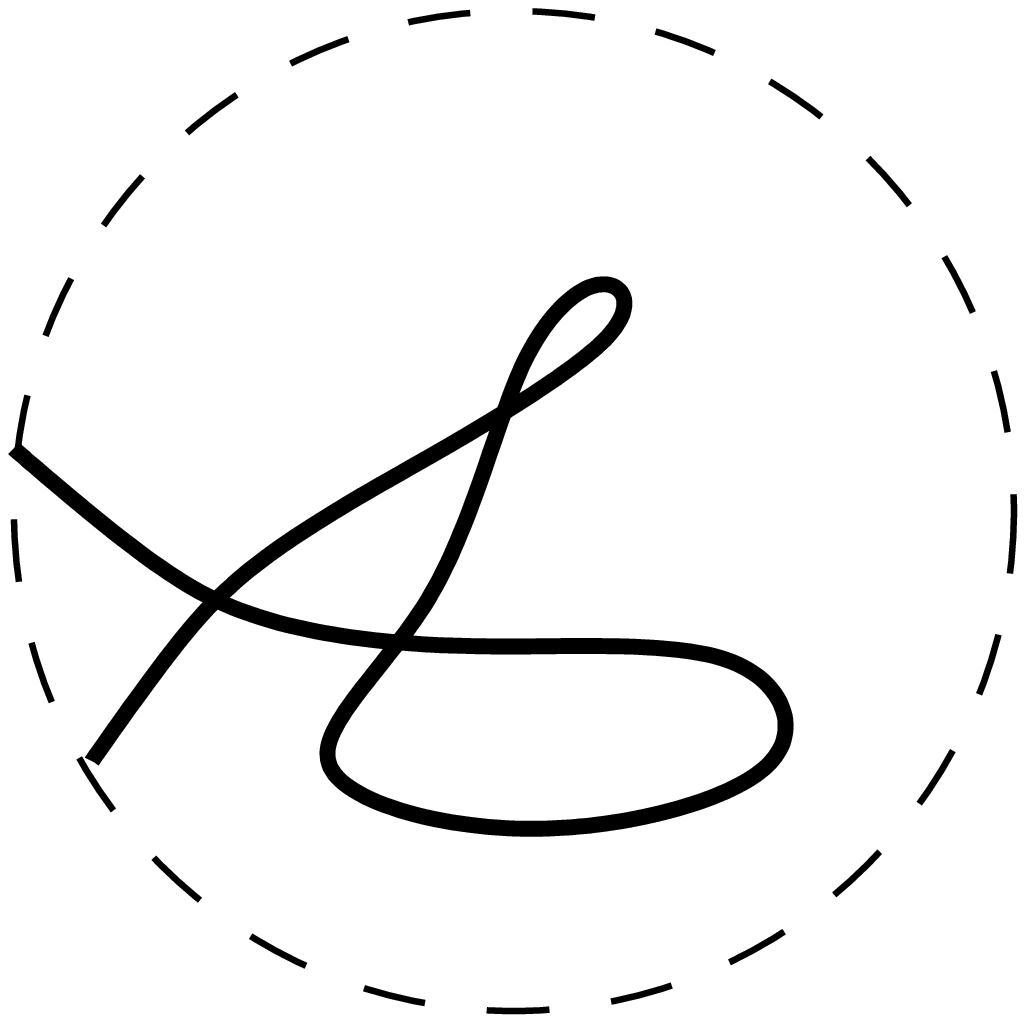}}
        \end{picture}}
$$

\sbp{1.3. Definition} (\cite{AC-1, AC-2}). Let $x$ be
the horizontal coordinate on the disk $B$ and $y$ be the vertical
coordinate. {\it A divide link} $\L_D$ is a link in the 3-sphere
$S^3=\{(x,y,u,v)\in\R^4 | x^2+y^2+u^2+v^2=1\}$ such that $(x,y)$ is a
point on $D$ and $u$, $v$ are the coordinates of a tangent vector to
$D$ at the point $(x,y)$.

So each interior point of $D$ has two corresponding points on $\L_D$,
whereas a boundary point of $D$ gives a single point on $\L_D$.

The number of components of $\L_D$ equals to the number of branches of
the divide $D$ which is the number of the copies of the unit interval
$I$ in definition {\bf 1.1}. In particular, if $D$ consists of only one
branch (like in {\bf 1.2}) then $\L_D$ will be a knot.

\sbp{1.4. Properties of $\L_D$.}

\sbp{1.4.1.} Topological type of a divide link does not change under
a regular transversal isotopy of the disk $B$. So it does not
depend on the choice of coordinates in 1.3. Also it does
not change under a moving of a piece of the curve $D$ through a
triple point \cite{CP}. In particular, the following two divides
have the same knot type as in 1.2.
$$
  \mbox{\begin{picture}(300,60)(0,0)
  \put(0,0){\epsfxsize=60pt \epsfbox{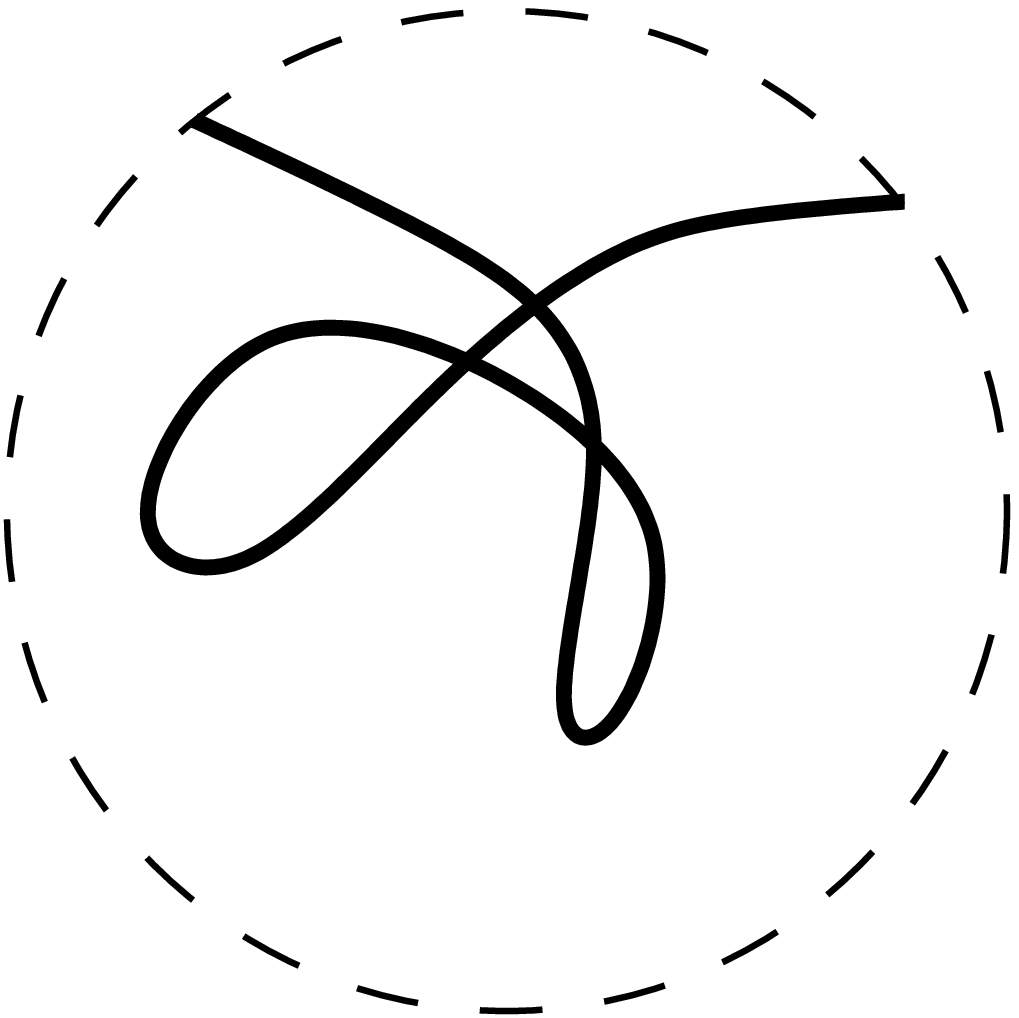}}
  \put(200,0){\epsfxsize=60pt \epsfbox{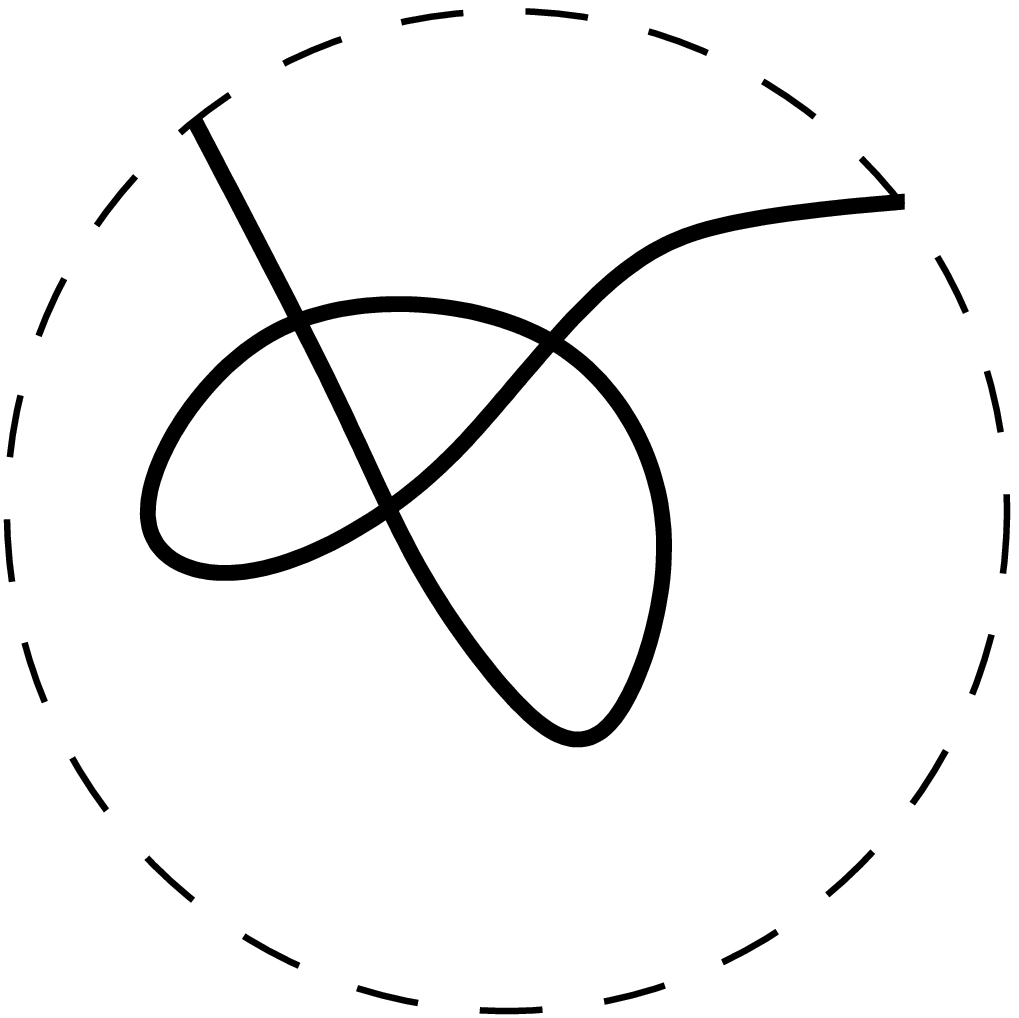}}
        \end{picture}}
$$

\sbp{1.4.2.} The link $\L_D$ has a natural orientation. Indeed, choose any
orientation of every branch of $D$. Let $(u,v)$ be the tangent vector
to $D$ at $(x,y)$ pointing to the direction of the chosen orientation
of $D$.  Then the orientation of $\L_D$ is given by the vector
$(\dot{x},\dot{y},\dot{u},\dot{v})$. It is easy to see that this
orientation of $\L_D$ does not depend on the choice of orientations of
branches of $D$.

\noindent\parbox[t]{370pt}{ \sbp{\quad 1.4.3.} In \cite{AC-1}
A'Campo proved that any singularity link is a divide link. But
the divide in the picture on the right does not occur in singularity
theory. Indeed, by the real morsification method one can find
that the Milnor number of the corresponding singularity would be 4. However
there are only two singularities with Milnor number 4, $A_4$ and
$D_4$. But their intersection forms are different from the one
which corresponds to this divide. The knot arising from
divide is $10_{145}$ (see \cite{AC-2,Ch}). }\qquad\vspace{8pt}
\parbox[t]{80pt}{
$$
  \mbox{\begin{picture}(60,65)(0,0)
  \put(0,-10){\epsfxsize=60pt \epsfbox{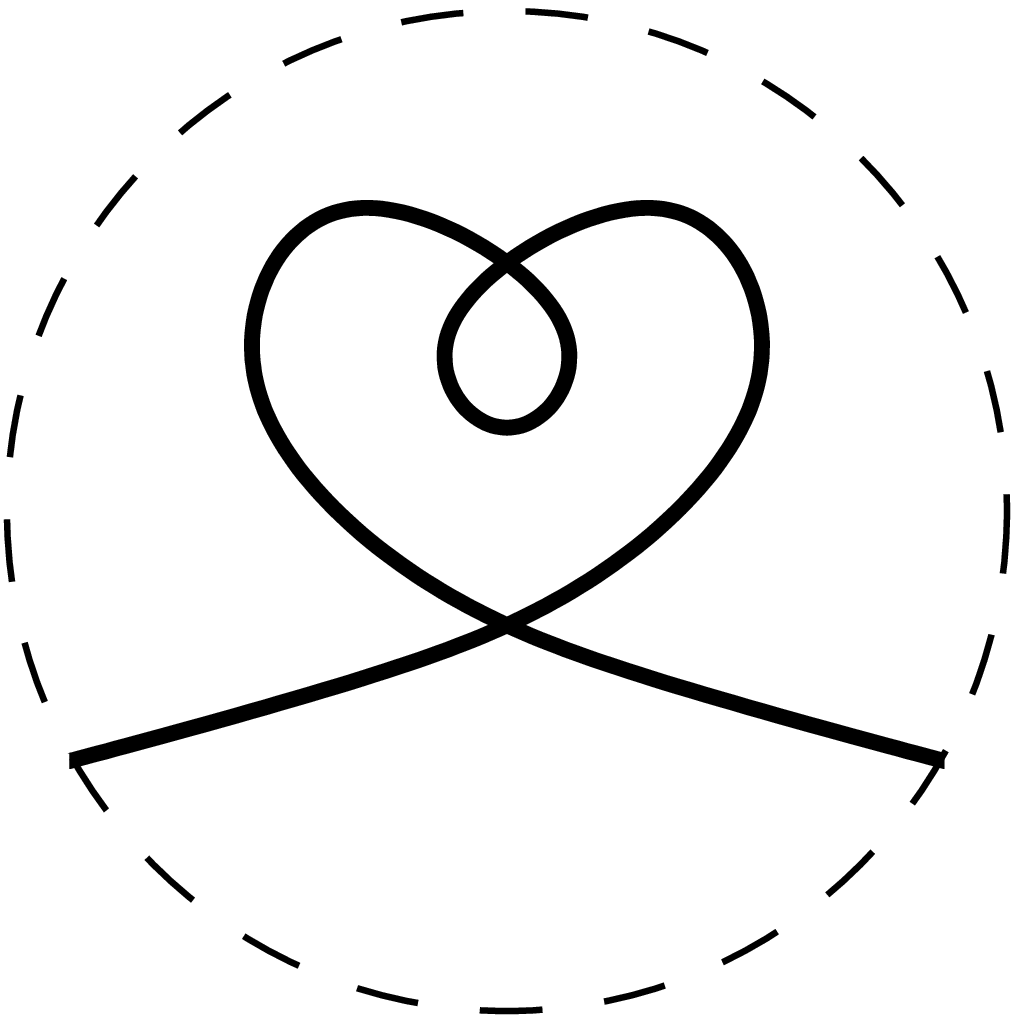}}
        \end{picture}}
$$}

It is known for a long time \cite{Bu} that all singularity knots
are classified by the Alexander polynomial.  N.A'Campo
(\cite{AC-3}) found two different divide knots with the same
Alexander polynomial. H.Morton \cite{Mor}, whom I have shown 
A'Campo's example, found that these knots are mutant. So they
cannot be distinguished by any classical polynomial invariants
(Jones, HOMFLY, Kauffman). He distinguished them by a quantum
invariant coming from the Lie algebra $\mathfrak{gl}_n$ 
in a certain
higher (non standard) representation.

\sbp{1.4.3.}
A'Campo \cite{AC-2} showed that the links $\L_D$
corresponding to a connected divide $D$ are fibered and computed
their monodromy in terms of combinatorics of the divide $D$.
Not all fibered links have the form $\L_D$. Figure eight knot $4_1$ is
not a divide knot. It is not clear how large is the class of
divide links in the class of all fibered links.

\sbp{1.4.4.} If we consider $\R^4$ with coordinates $x,y,u,v$ as a
complex plane $\C^2$ with coordinates $z_1=x+iu$ and $z_2=y+iv$,
then every tangent space to the unit sphere $S^3$ contains a
unique complex line which is a two-dimensional real subspace in the
tangent space. The distribution of these 2-planes forms the
standard contact structure on $S^3$. Divide knots are transversal
knots with respect to this contact structure. This fact was
noticed in \cite{AC-2} but the arguments given there should be modified.
For a one-branch divide $D$ the Bennequin number (self-linking
number) of the knot $\L_D$ is equal to $2\delta-1$, where $\delta$
is the number of double points of $D$.

\sbp{1.4.5.} N.A'Campo \cite{AC-2} proved that the unknotting
number of a one-branch divide knot $\L_D$ and the genus of $\L_D$
are equal to the number $\delta$ of double points of $D$. For a
two-branch divide $D$ the link $\L_D$ will have two components.
Their linking number is equal to the number of common double
points of the two branches of $D$.

\section{Casson's invariant of knots.}

\sbp{2.1. The Conway polynomial.} The subject of this section is well
known (see, for example \cite{Ka,PV}). The Conway polynomial $C(\L)$
of a link $\L$ a polynomial in a single variable $z$. It is one of the
simplest invariants of links in $\R^3$, defined by the two
properties: it is equal to one on the unknot and satisfies the skein relation
\def\sk#1{\mbox{\begin{picture}(30,10)(0,0)
        \put(0,-6){\epsfxsize=30pt \epsfbox{#1.eps}}
        \end{picture}}}
$$C(\ \sk{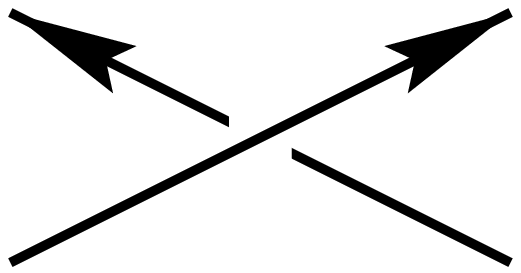}\ ) - C(\ \sk{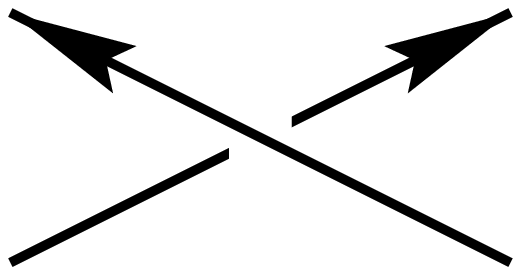}\ ) = z\cdot C(\ \sk{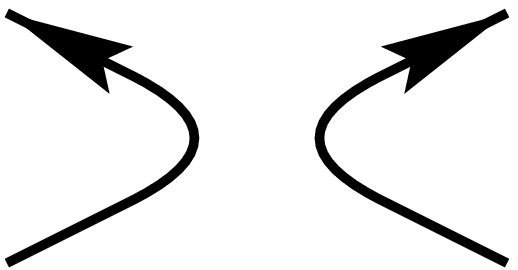}\ )\ ,
$$
where the three links are identical outside a small ball in $\R^3$
and look as shown inside the ball.

If $\L$ is a knot then the  Conway polynomial is even:
$$C(\L) = 1 + C_2(\L)\cdot z^2 + C_4(\L)\cdot z^4 + \dots
            + C_{2n}(\L)\cdot z^{2n}\ .$$

\sbp{2.2. Casson's invariant.} The coefficient $C_2(\L)$ is called
{\it Casson's knot invariant} (see \cite{PV}). It can be defined
(\cite{Ka}) by the
initial condition $C_2(unknot)=0$ and the skein relation:
$$C_2(\ \sk{sk1}\ ) - C_2(\ \sk{sk2}\ ) = lk(\ \sk{sk3}\ )\ ,$$
where $lk$ means the linking number of the two component link in the right
hand side of the relation.

The $mod\ 2$ reduction of Casson's invariant $C_2(\L)$ is called
{\it Arf invariant} of a knot $\L$.


\section{Arnold's invariants of immersed plane curves.}

\sbp{3.1. Arnold's invariants of divides.} In \cite{Ar} V.I.Arnold
defined three basic invariants $J^+$, $J^-$, $St$ of a generic
{\it closed} immersed plane curve. See a survey of various
explicit formulas for these invariants in \cite{CD}. Following
\cite{GZN} we can define these invariants for a one-branch divide
$D$ as the arithmetic mean of the corresponding Arnold's
invariants on the two closed curves obtained by smoothing the
union of $D$ with either arc of the unit circle. 
Invariants of divides thus defined will be denoted by the same symbols.

\sbp{3.2. Example.} Let us compute $J^-$ on the divide of the example
{\bf 1.5.2}.
\def\jpi#1{J^-\left(\ \mbox{\begin{picture}(30,10)(0,0)
        \put(0,-9){\epsfxsize=30pt \epsfbox{#1.eps}}
        \end{picture}}\ \right)}
\def\jdpi#1{J^-\left(\ \mbox{\begin{picture}(30,10)(0,0)
        \put(0,-13){\epsfxsize=30pt \epsfbox{#1.eps}}
        \end{picture}}\ \right)}
$$\jdpi{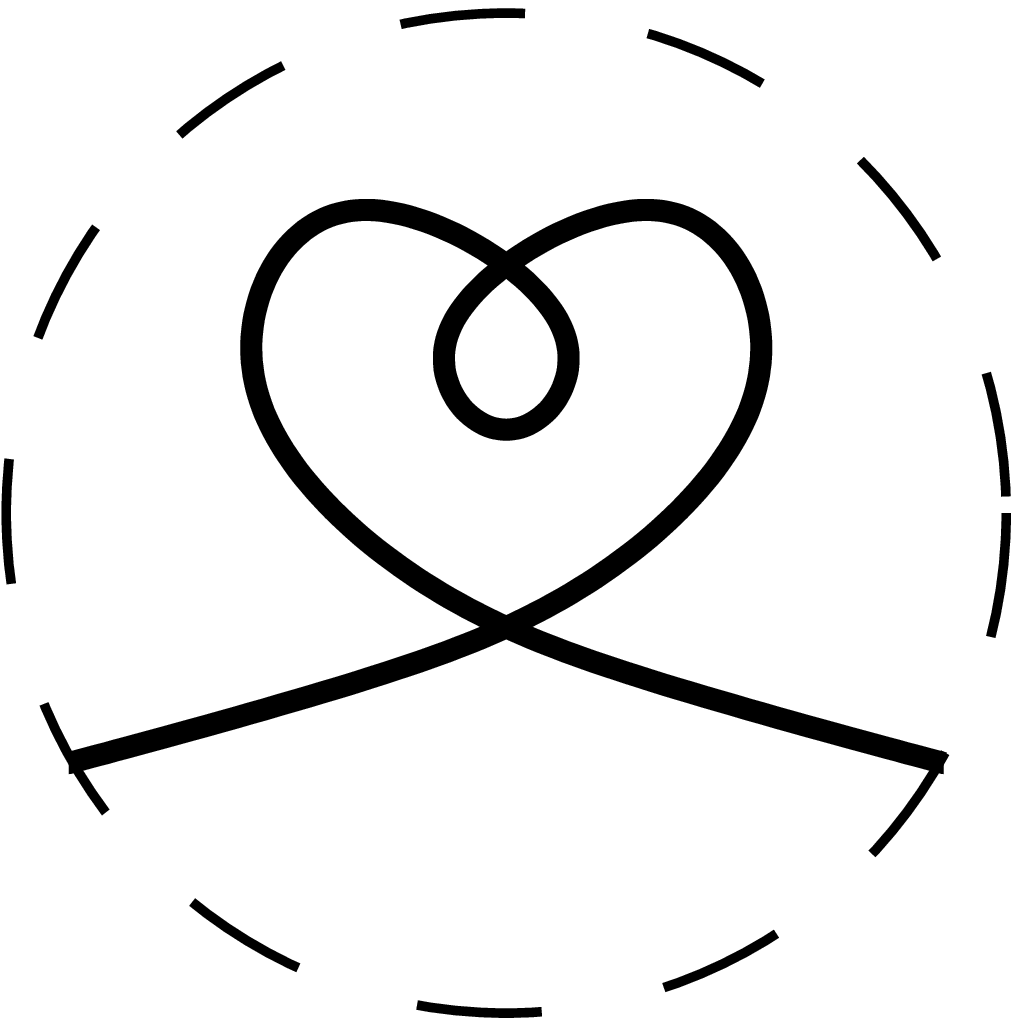}
 = {\raisebox{10pt}{$\jpi{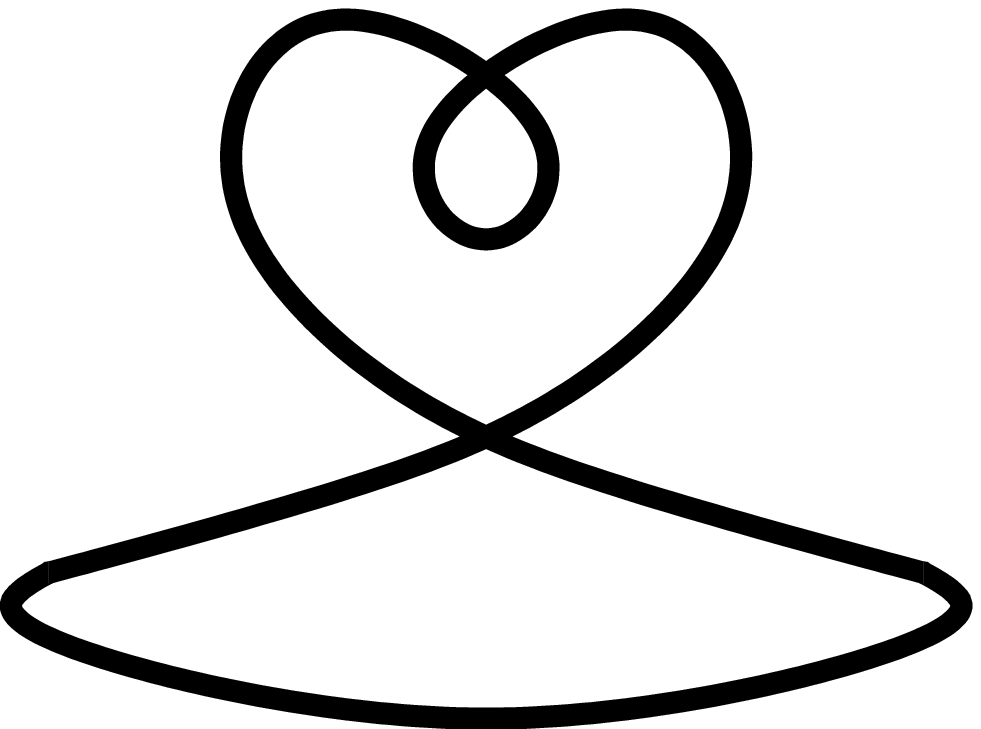} + \jpi{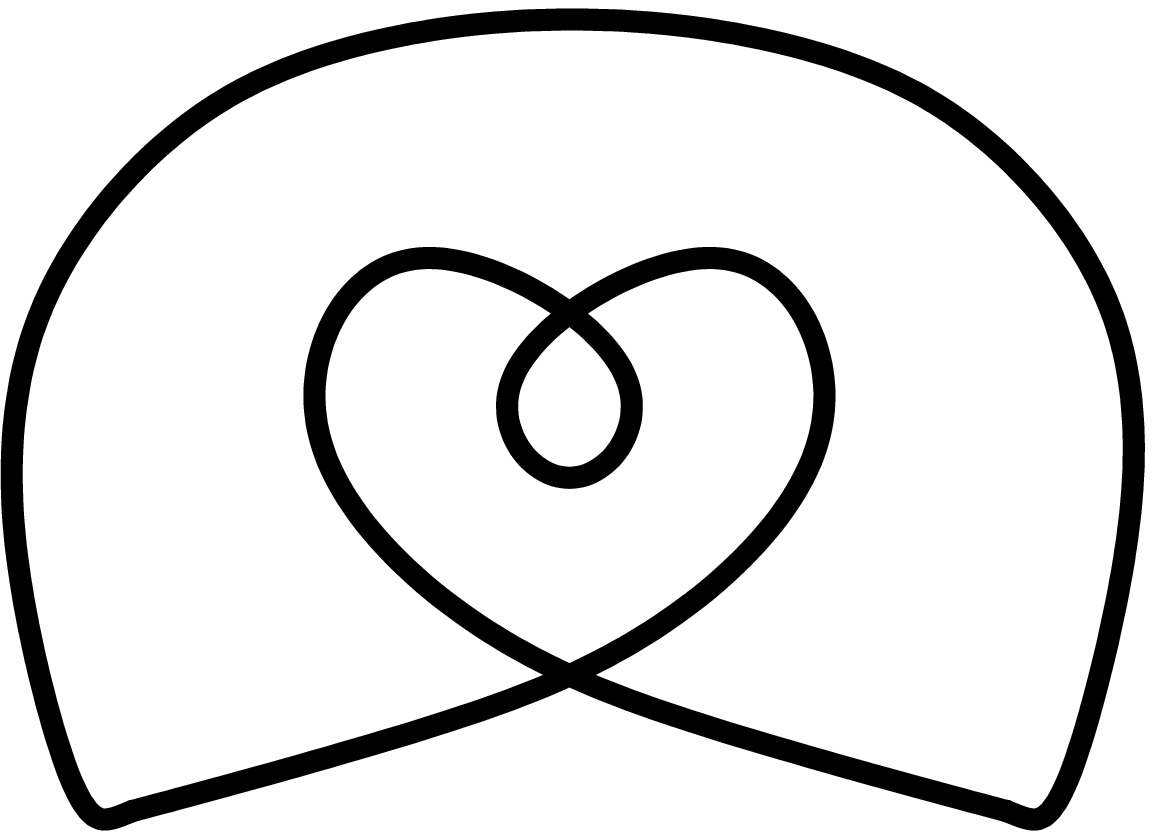}$} \over 2}\ .$$
According to Arnold's table \cite[page 14]{Ar}, the first term of the right
hand side equals $-4$ and the second term equals $-8$. So the result 
is $(-4-8)/2=-6$.

\sbp{3.3. $J^\pm$ as invariants of order 1.}
We are going to use a description of Arnold's $J^\pm$ invariants
as invariants of order 1 in the sense of the theory of finite type
invariants.

Let us fix the images of the ends points of the interval $I$  at
the boundary of the disk $\partial B$. A version of the classical
theorem of H.~Whitney states that the space of all smooth
immersions of the interval I in the disk $B$ mapping $\partial I$
into two fixed points is made up of a countable number of
connected components differing by the absolute value of the rotation
number.

Choose a standard divide $D_i$ for each nonnegative value of the rotation
number $i$\vspace{10pt}:
\def\stdi#1{\begin{array}{c}
             \mbox{\begin{picture}(30,10)(0,0)
             \put(0,-13){\epsfxsize=30pt \epsfbox{st#1.eps}}
             \end{picture}}\vspace{10pt} \\ D_{#1} \end{array}}
$$\stdi{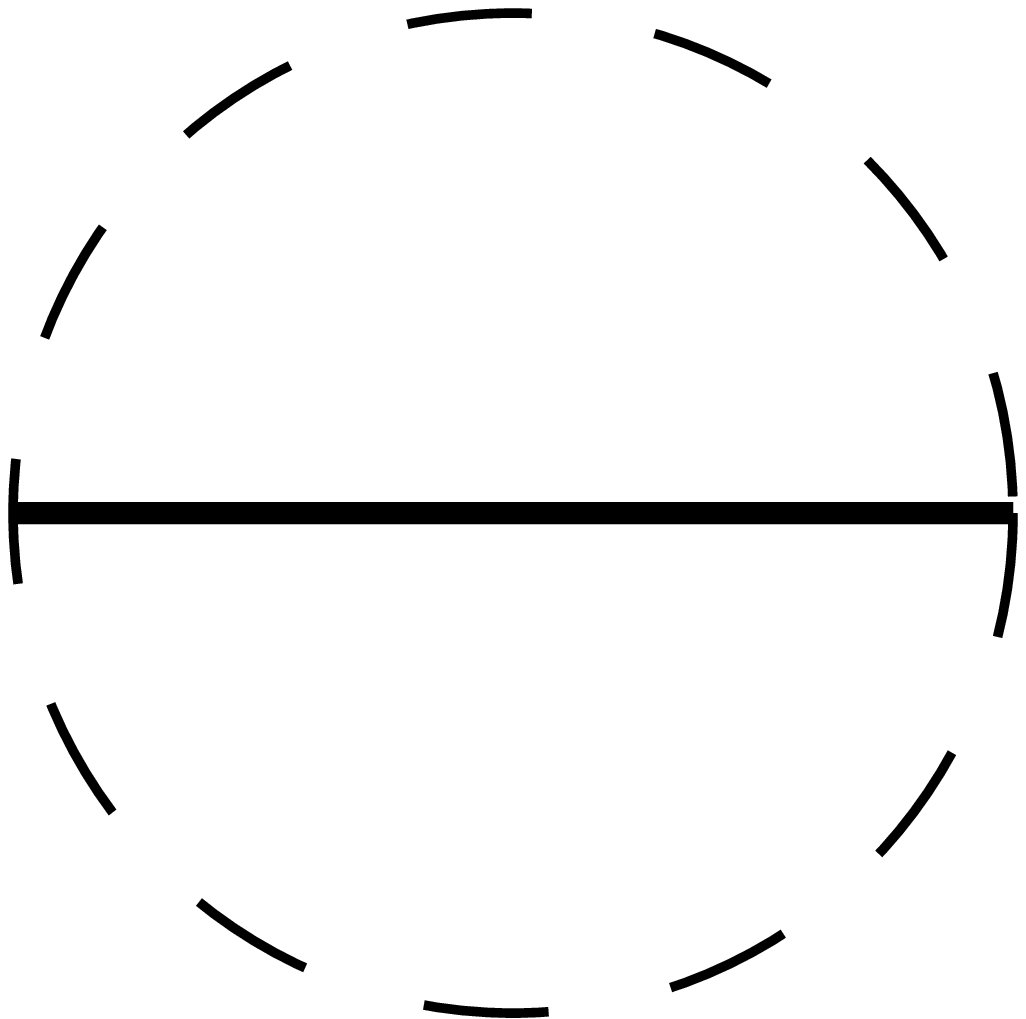}\qquad \stdi{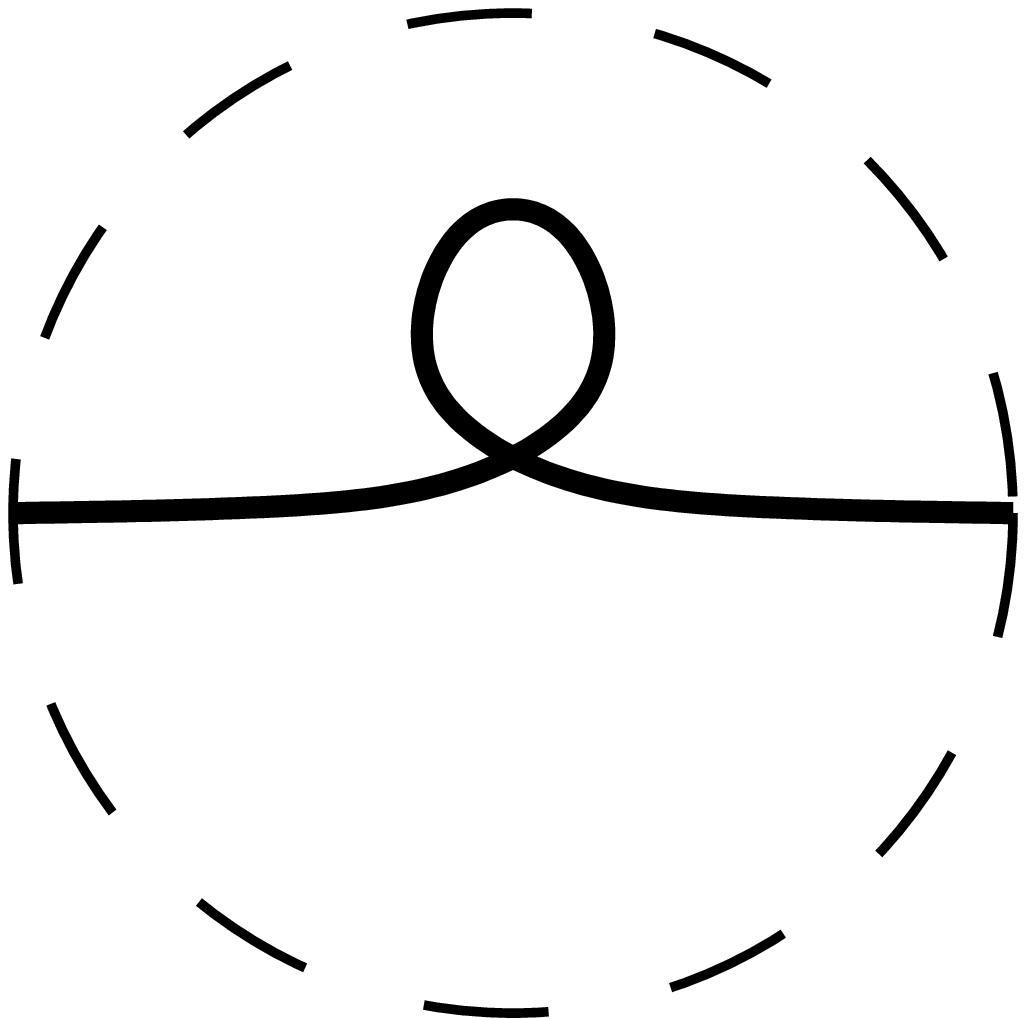}\qquad \stdi{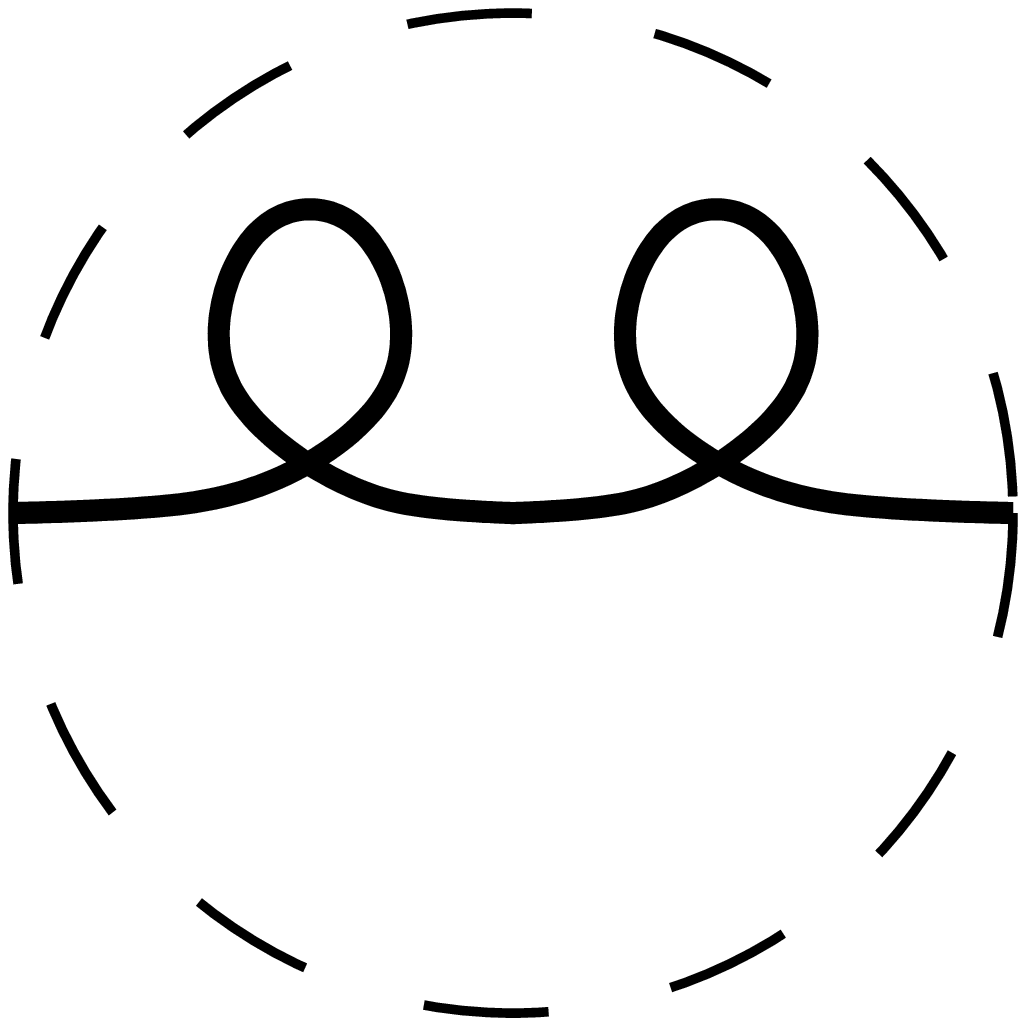}\qquad \stdi{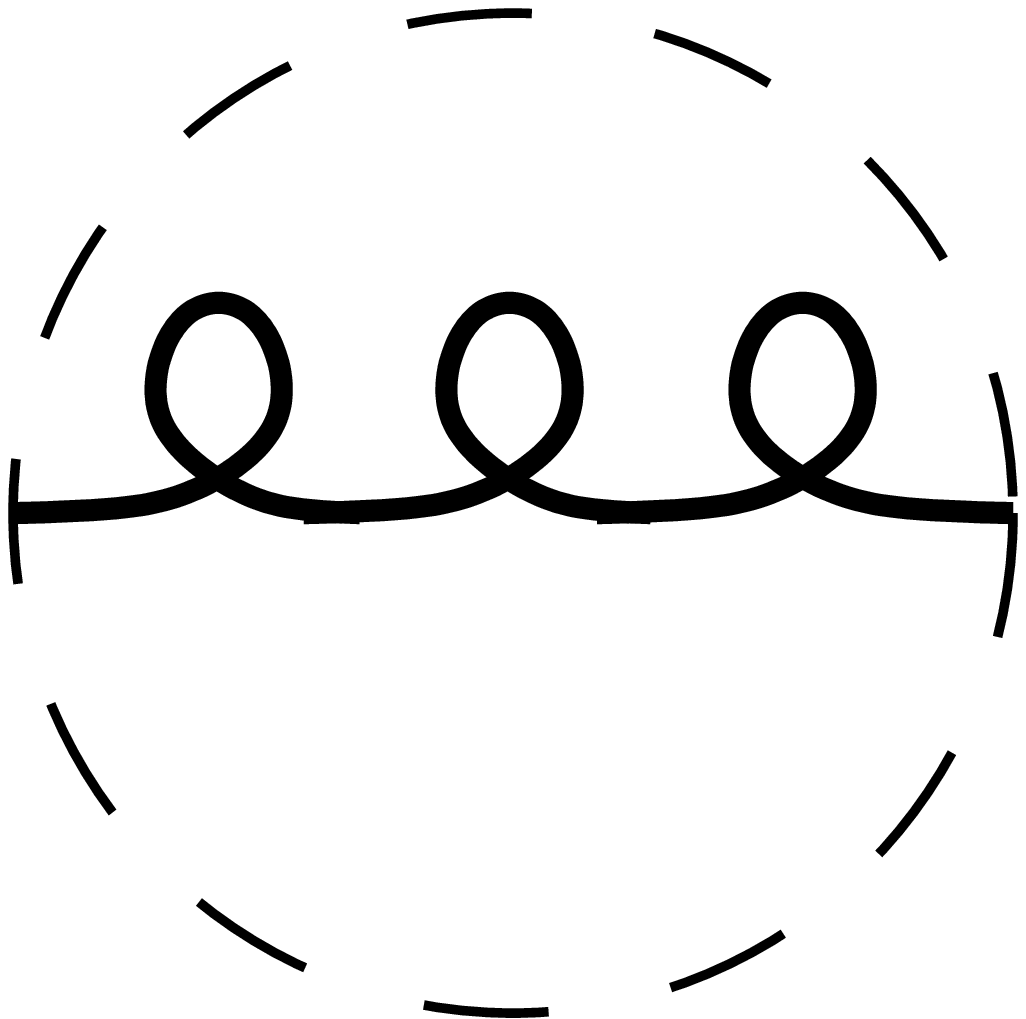}\qquad
\stdi{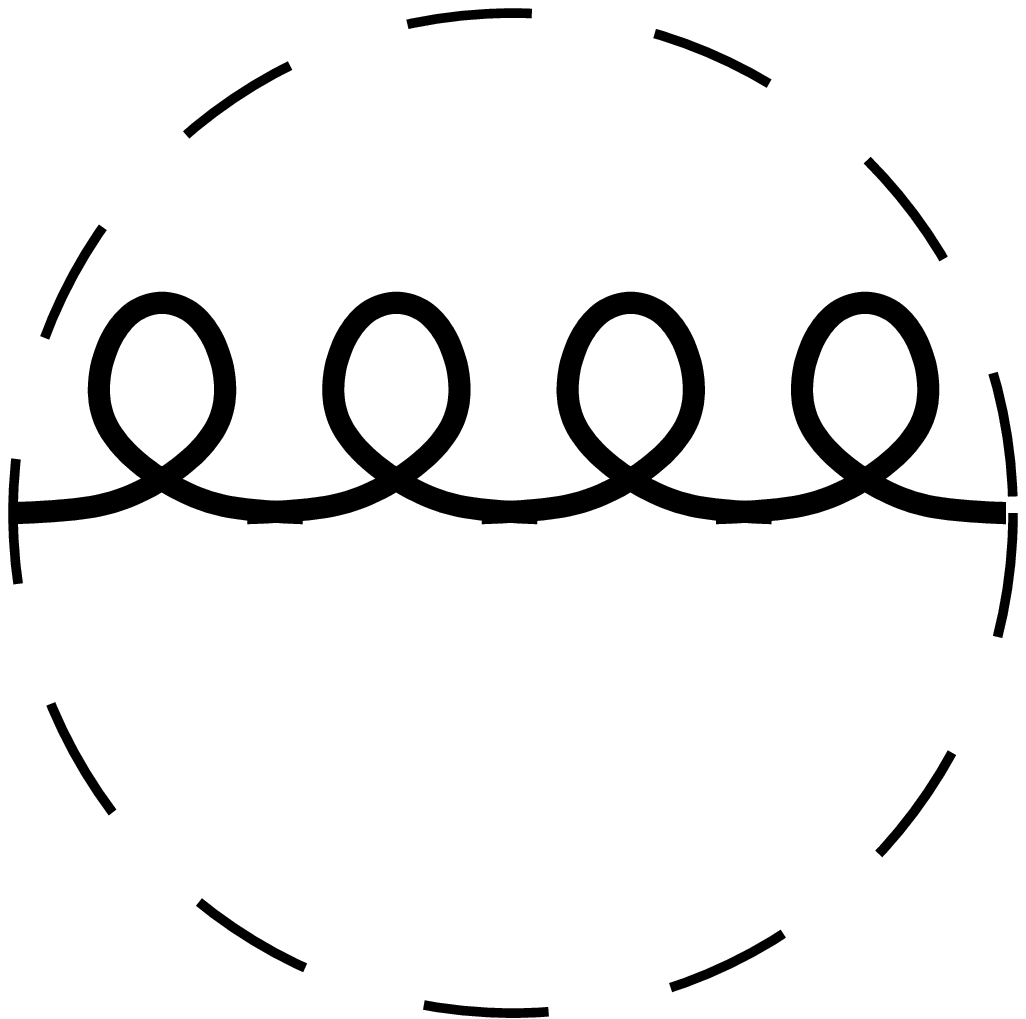}\qquad \stdi{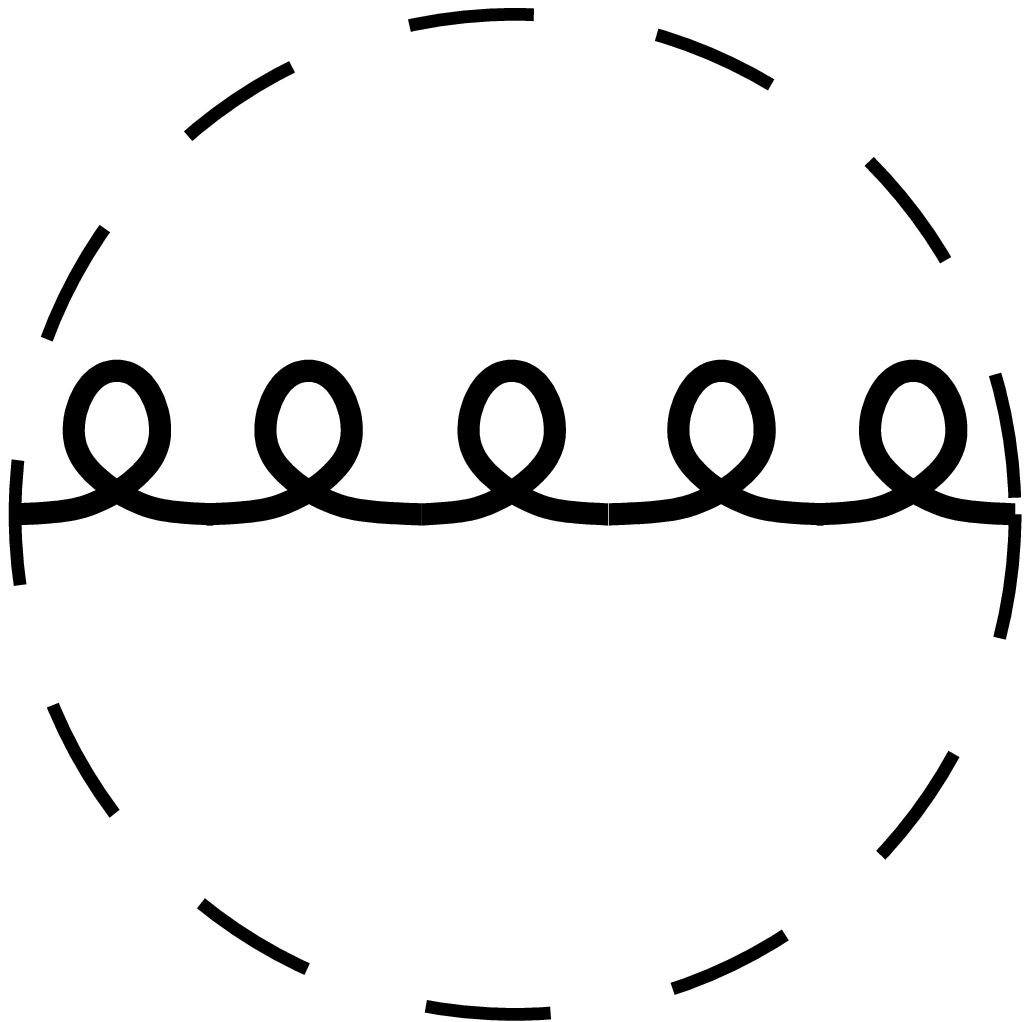} \qquad \raisebox{10pt}{\dots}$$ 
So every
one-branch divide can be deformed to one of $D_i$'s in the
space of immersions.

During the deformation at certain moments non-divides may occur
because of non allowed singularities of the corresponding curves.
For a generic deformation two types of such non-divides can occur:
either a {\it triple point} on the curve or a {\it self-tangency}
of the curve. In fact the self-tangency event can be split in two
types: {\it direct self-tangency}, where the two tangent strings
have coherent orientations (for any of the two possible
orientations of the curve), and {\it inverse self-tangency}, where
the two tangent strings have opposite orientations. A first order
invariant is defined by 
its jumps on the events of the three types above.
and 
its values on the standard divides 

In particular, $J^-$ and $J^+$ are defined by the following relations:
\def\jsk#1{J^-\left(\ \mbox{\begin{picture}(30,10)(0,0)
        \put(0,-6){\epsfxsize=30pt \epsfbox{#1.eps}}
        \end{picture}}\ \right)}
\def\jsp#1{J^+\left(\ \mbox{\begin{picture}(30,10)(0,0)
        \put(0,-6){\epsfxsize=30pt \epsfbox{#1.eps}}
        \end{picture}}\ \right)}
$$\begin{array}{ll@{\quad}|@{\quad}ll}
\mbox{(i)}^- & \jsk{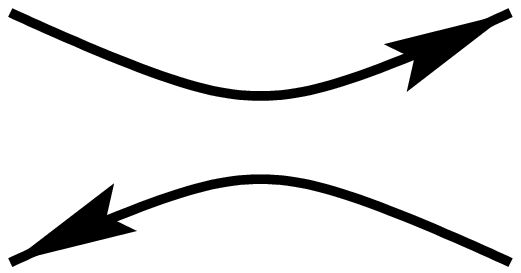} -\jsk{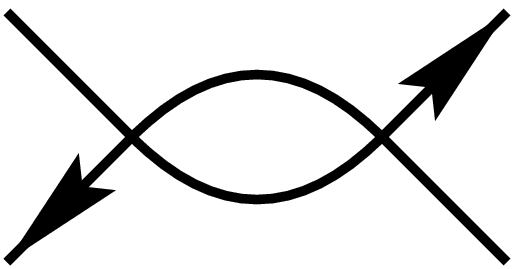} = 2\ ; &
\mbox{(i)}^- & \jsp{ib} =\jsp{ia}\ ; 
               \raisebox{-15pt}{\makebox(0,20){}}\\
\mbox{(i)}^+ & \jsk{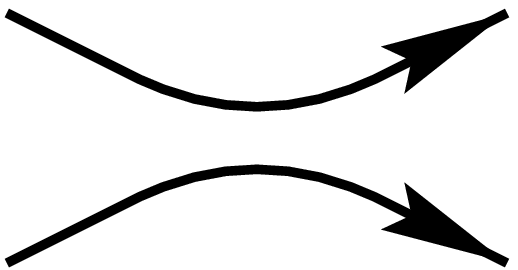} =\jsk{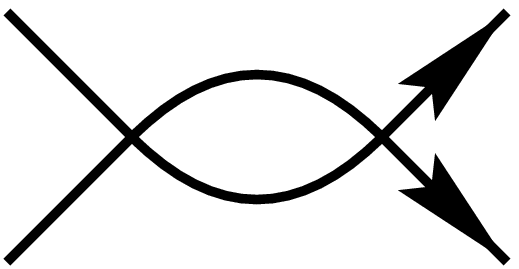}\ ; &
\mbox{(i)}^+ & \jsp{db} -\jsp{da} = -2\ ; 
               \raisebox{-15pt}{\makebox(0,20){}}\\
\mbox{(i)}^* & \jsk{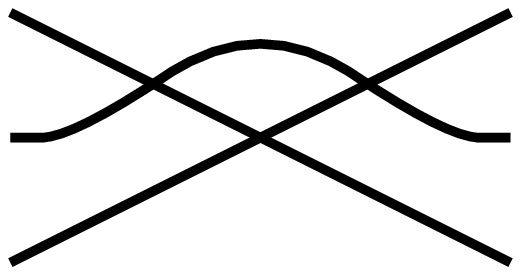} =\jsk{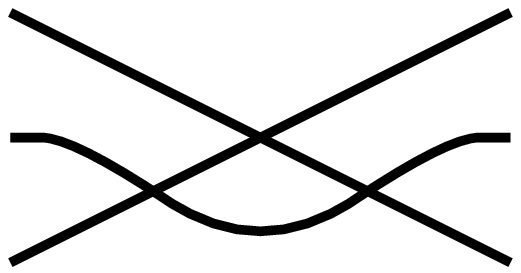}\ ; &
\mbox{(i)}^* & \jsp{tb} =\jsp{ta}\ ; 
               \raisebox{-15pt}{\makebox(0,20){}}\\
\mbox{(ii)} & J^-(D_i) = -2i\ ; &
\mbox{(ii)} & J^+(D_i) = -i\ ;  
\end{array}$$
Here in each equality we mean two divides identical outside a
small fragments which are explicitly shown.

\sbp{3.4. Actuality tables for $J^\pm$.}
According to V.Vassiliev \cite{Va} actuality tables provide a way to 
organize the data necessary for the computation of a single finite 
order invariant.

For the first order invariants $J^\pm$ the actuality tables are 
essentially the same thing as the set of equations (i)$^\bullet$--(ii). To 
describe the top row of the actuality table we use chord diagrams. 
Since in the $J^\pm$ theory we have two types of singular events, the 
inverse and direct self-tangencies, we need two types of chords. We 
use a dashed chord to depict the inverse self-tangency, and a solid 
chord to depict the direct self-tangency.
\def\jd#1{\ \mbox{\begin{picture}(30,10)(0,0)
        \put(0,-6){\epsfxsize=30pt \epsfbox{#1.eps}}
        \end{picture}}\ }
\def\jcd#1{\ \mbox{\begin{picture}(40,10)(0,0)
        \put(0,-5){\epsfxsize=40pt \epsfbox{#1.eps}}
        \end{picture}}\ }
$$\jcd{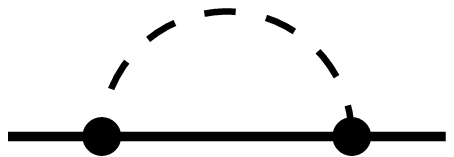}=\jd{ib}-\jd{ia}\qquad\qquad
\jcd{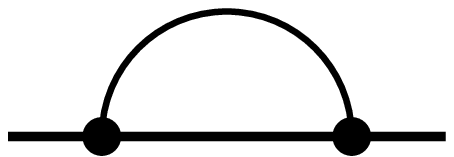}=\jd{db}-\jd{da}$$

The actuality tables for $J^\pm$ look as follows
$$\begin{array}{|ll@{\qquad}|@{\qquad}ll|} 
\hline
\mbox{(i)}^- & J^-\Big(\jcd{m}\Big)= 2 &
\mbox{(i)}^+ & J^+\Big(\jcd{p}\Big)= -2 
       \raisebox{-12pt}{\makebox(0,30){}} \\ \hline
\mbox{(ii)} & J^-(D_i)= -2i &
\mbox{(ii)} & J^+(D_i)= -i 
       \raisebox{-8pt}{\makebox(0,22){}} \\ \hline
\end{array}$$

\sbp{3.5. Example.} Let us compute $J^-$ again for the example
{\bf 1.5.2} but using now the relations (i)--(iv) only. We introduce 
the orientation as shown in order to distinguish between direct 
(i)$^+$ and inverse (i)$^-$ self-tangencies.
\def\en#1{\raisebox{-5pt}{$\begin{array}{cc} = \vspace{-5pt}\\
                             \mbox{\scriptsize #1}\end{array}$}}
$$\begin{array}{ccl}
\jdpi{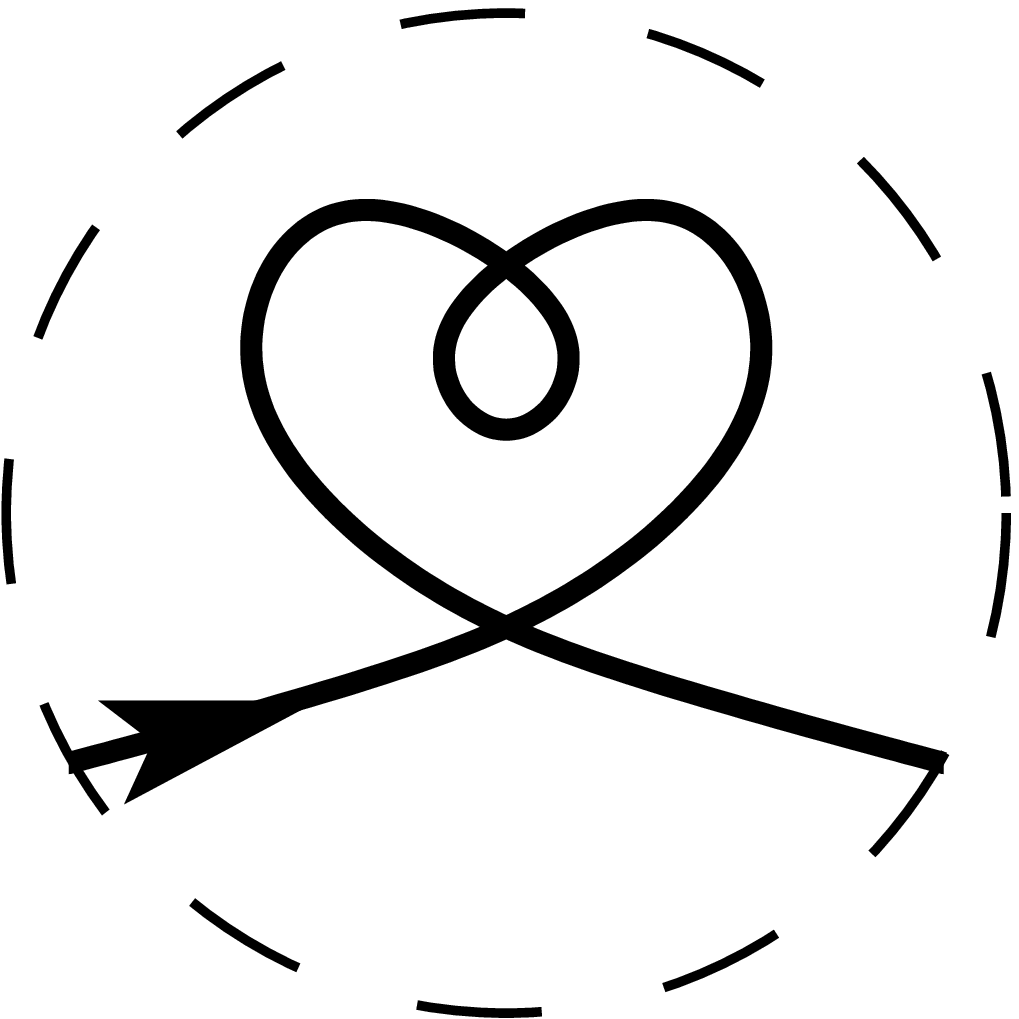} &=& \jdpi{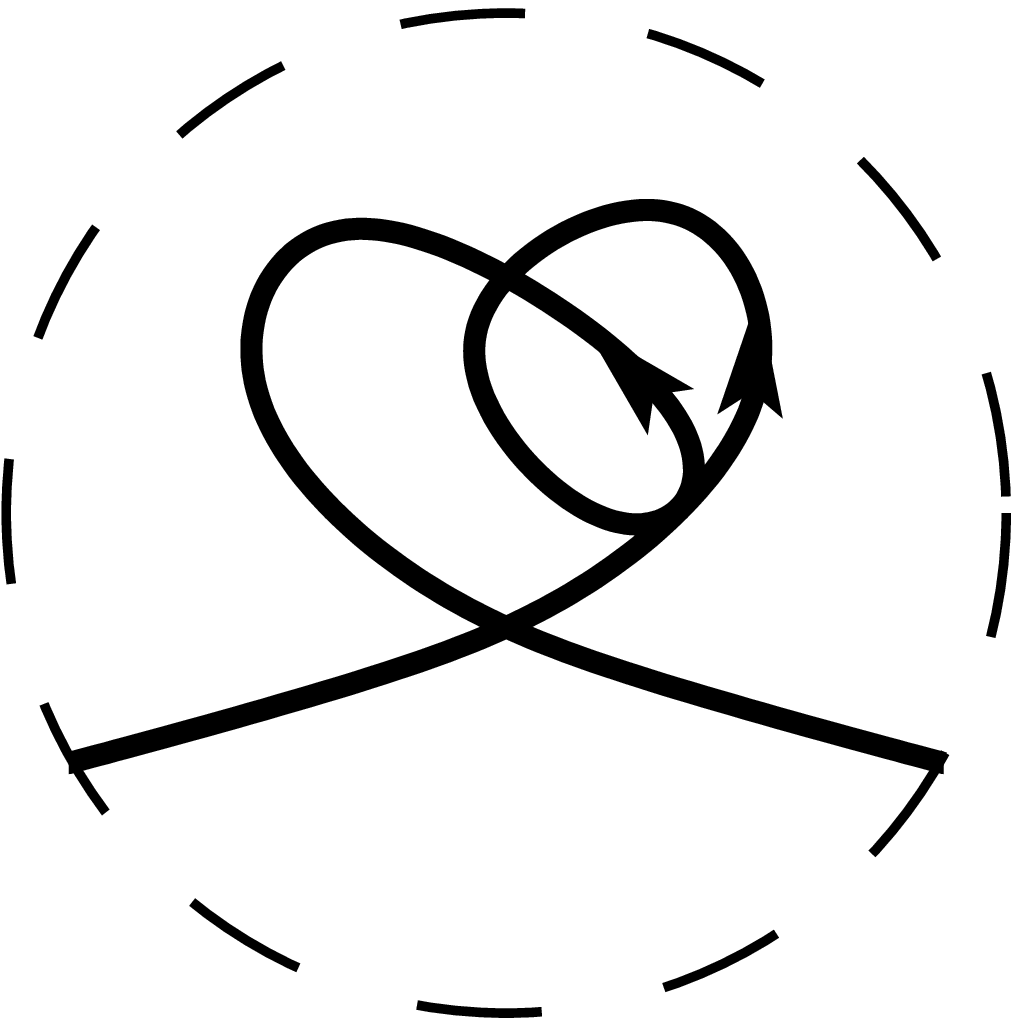} + \jdpi{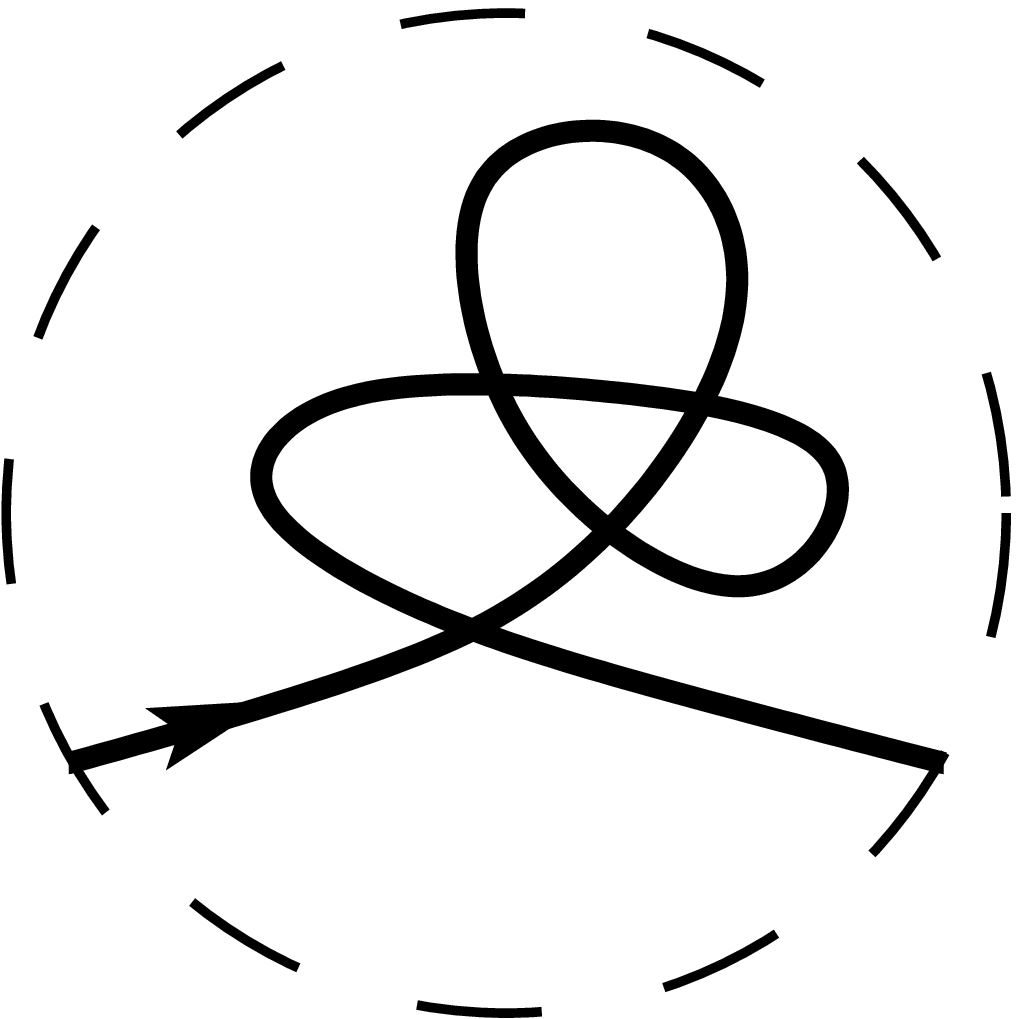} \en{(i)$^+$}
                     \jdpi{ex35-2} \en{(i)$^*$}
                  \jdpi{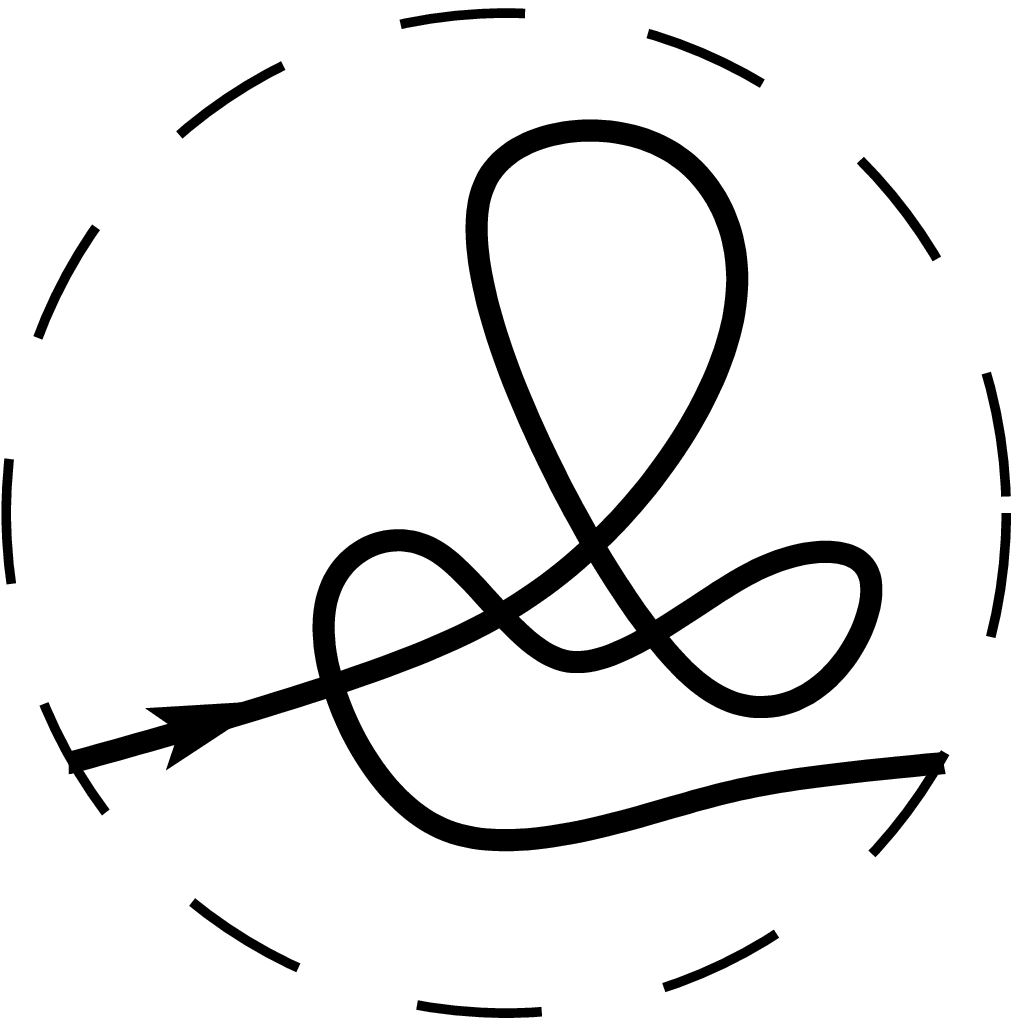} \vspace{20pt}\\
        &=& -\jdpi{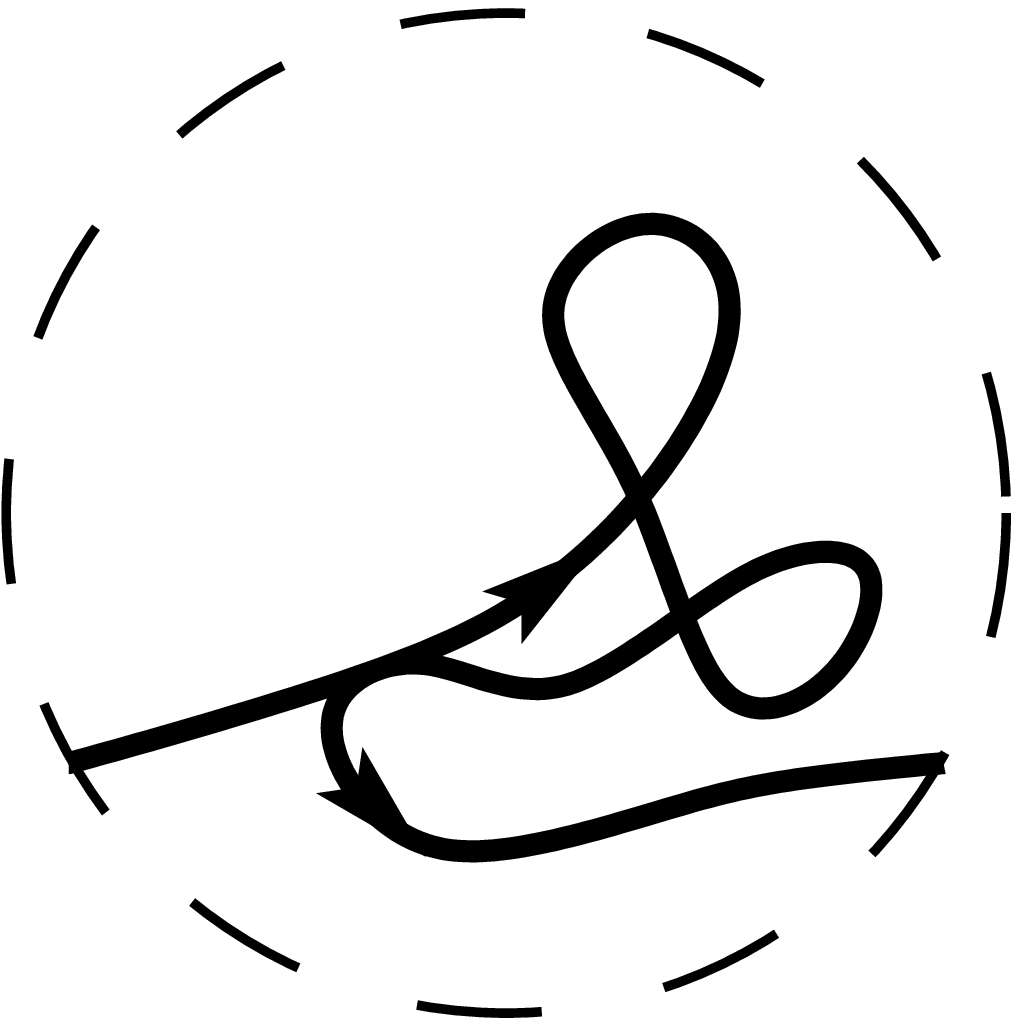} + \jdpi{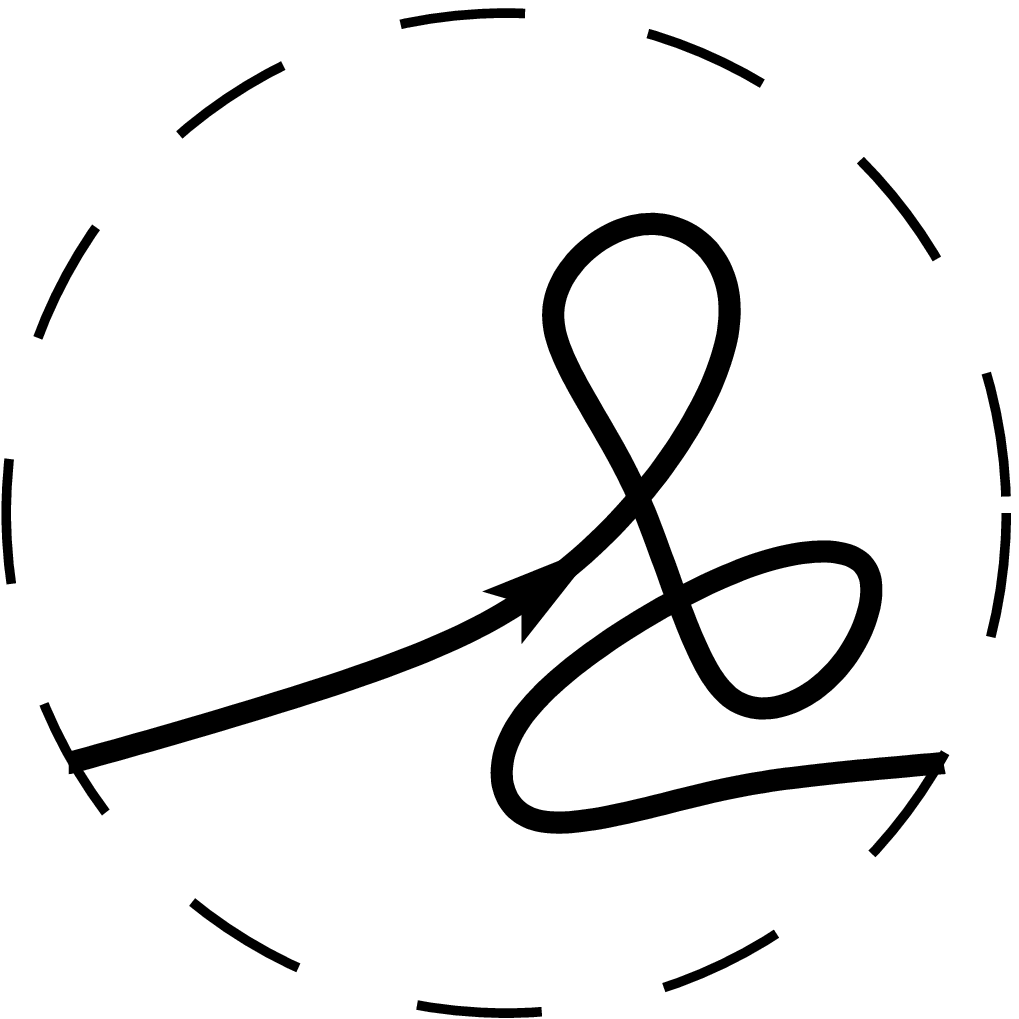} \en{(i)$^-$} -2+\jdpi{st2} 
                      \vspace{20pt}\\
              &=& -2+J^-(D_2) \en{(ii)} -2-2\cdot2=-6\ .
\end{array}
$$

\section{Main result}

\sbp{4.1. Theorem.} {\it For a one-branch divide $D$ the Casson 
invariant $C_2(\L_D)$ is equal to the invariant $J^\pm_2(D)$ defined 
below in {\bf 4.2}.}

\sbp{Corollary} (Gusein-Zade--Natanzon \cite{GZN}). 
{\it For a one-branch divide $D$, 
$\Arf(\L_D) = J^-(D)/2 (mod\ 2)$.}

\sbp{4.2. Definition of the invariant $J^\pm_2$.}\makebox(0,0){}

The invariant $J^\pm_2$ is a second order $J^\pm$-type
invariant. In particular, it does not change under the triple point 
move. To define
it we need to define its values on the chord diagrams with two
chords and also to define its values on canonical divides with at
most one self-tangency point. We have chosen the canonical divides
without self-tangencies in {\bf 3.3}. Now we must choose the
canonical divides with one self-tangency in such a way that any
divide with a single self-tangency can be deformed to the
canonical one if we allow to pass through the codimension two strata 
in the space of immersions.

\sbp{4.2.1. Canonical divides with a self-tangency point.} Our
choice of the canonical divides is:
$$D^-_{m,n} =
\begin{picture}(180,20)(0,0)
  \put(0,-10){\epsfxsize=150pt \epsfbox{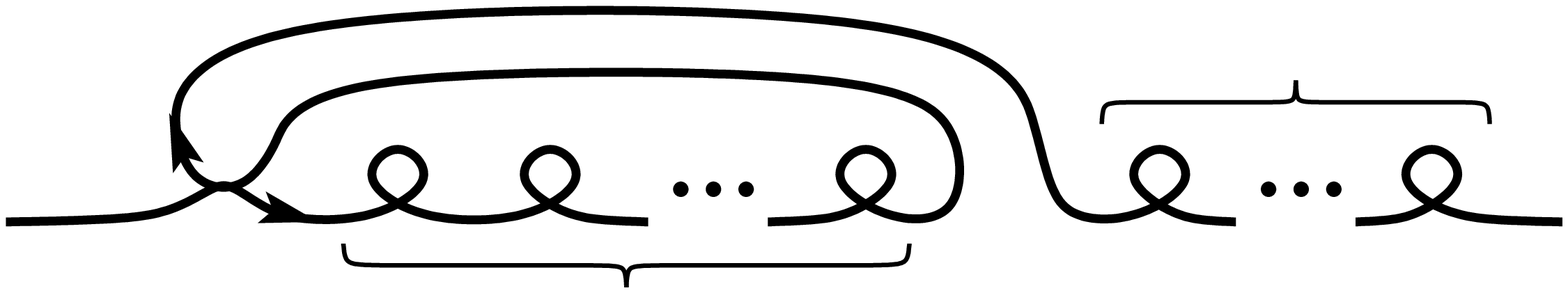}}
  \put(45,-15){\mbox{\scriptsize $m$ curls}}
  \put(110,12){\mbox{\scriptsize $n$ curls}}
\end{picture}
\qquad D^+_{m,n} =
\begin{picture}(180,20)(0,0)
  \put(0,-10){\epsfxsize=150pt \epsfbox{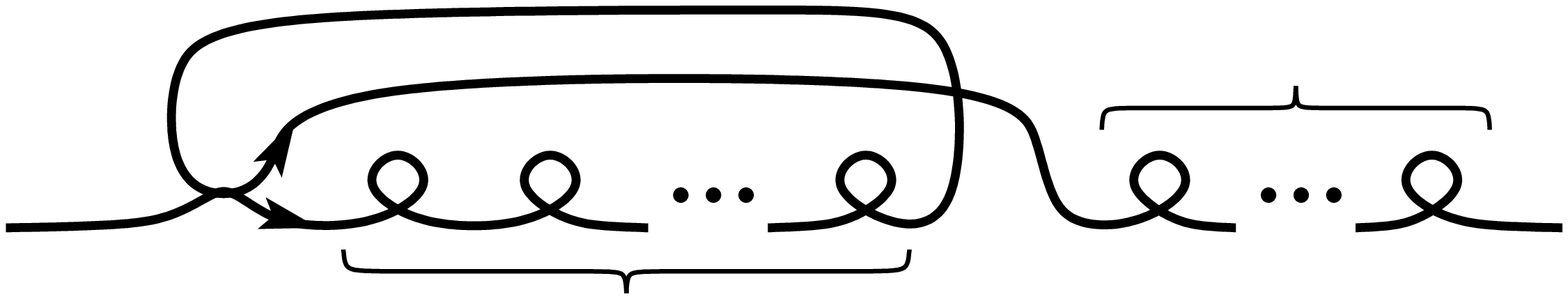}}
  \put(45,-15){\mbox{\scriptsize $m$ curls}}
  \put(110,12){\mbox{\scriptsize $n$ curls}}
\end{picture}
\vspace{8pt}$$ Here we omit the boundary circle (early 
depicted as a dashed circle) of the unit disk containing the divide.
The parameters $m$ and $n$ in these divides run over all integral
numbers. Under a negative curl we mean the curl which is going
clockwise (instead of counterclockwise as above). Here is a couple
of examples.
$$D^-_{-3,2} =
\begin{picture}(150,15)(0,0)
  \put(0,-10){\epsfxsize=120pt \epsfbox{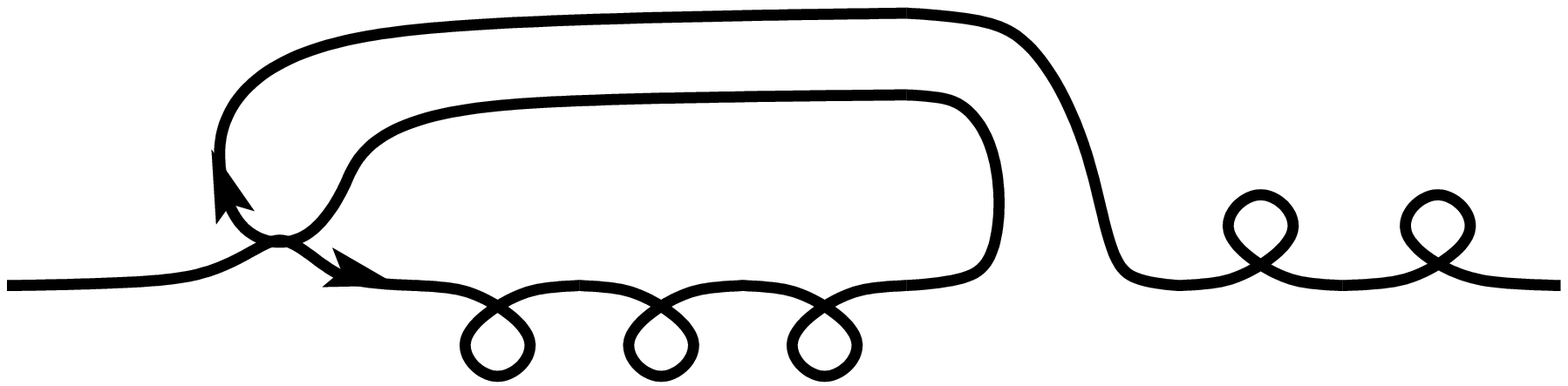}}
\end{picture}
\qquad D^+_{2,-3} =
\begin{picture}(150,15)(0,0)
  \put(0,-10){\epsfxsize=120pt \epsfbox{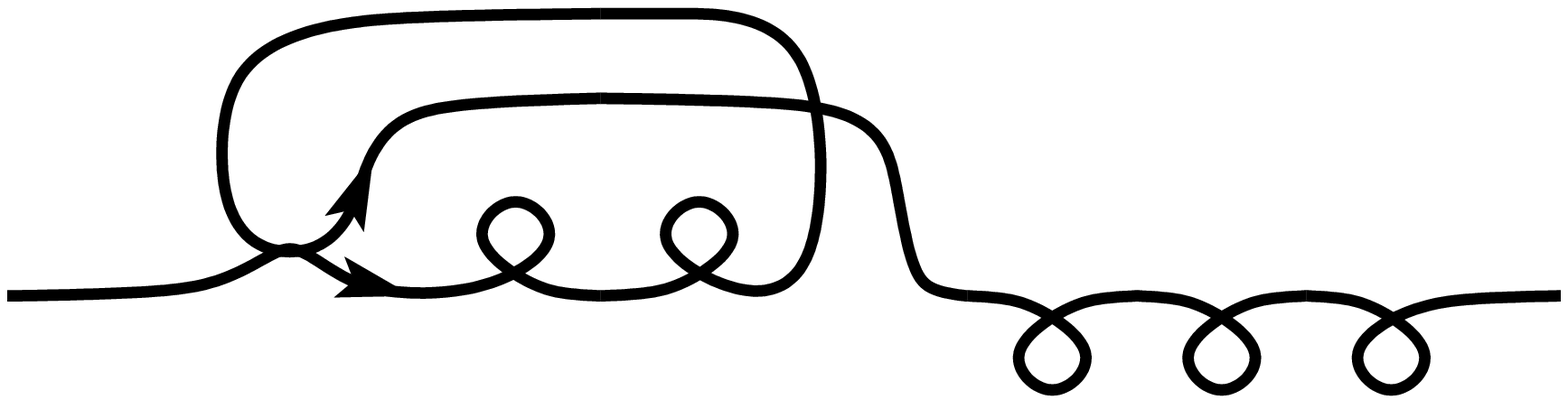}}
\end{picture}
\vspace{8pt}$$

\sbp{4.2.2. Actuality table for $J^\pm_2$.} The actuality table for
$J^\pm_2$ looks as follows.
\def\jcd#1#2{J^\pm_2\Big(\ \begin{picture}(60,20)(0,0)
  \put(0,0){\epsfxsize=60pt \epsfbox{#1.eps}}
\end{picture}\ \Big) = #2}
$$\begin{array}{|c@{\qquad}c@{\qquad}c|}
\hline&&\vspace{-8pt}\\
\jcd{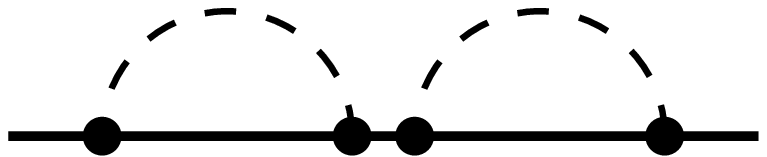}{0} & \jcd{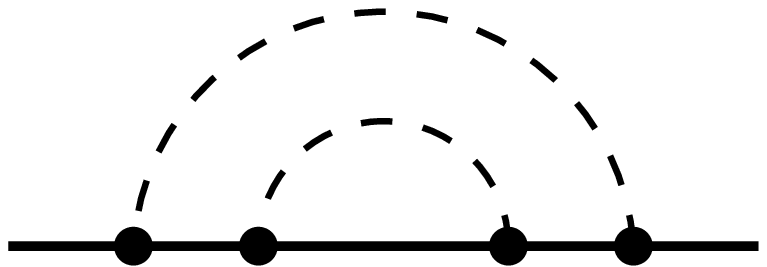}{4} & \jcd{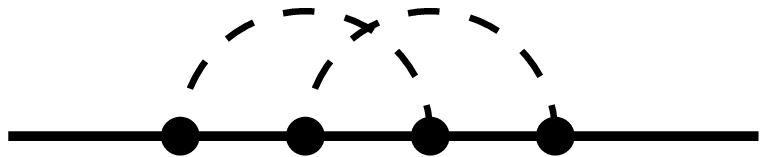}{2} \\
\jcd{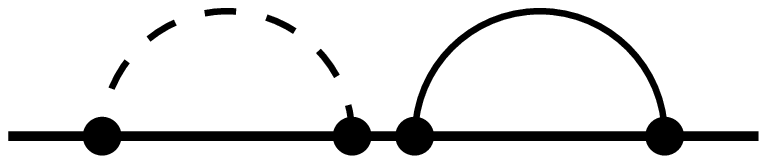}{0} & \jcd{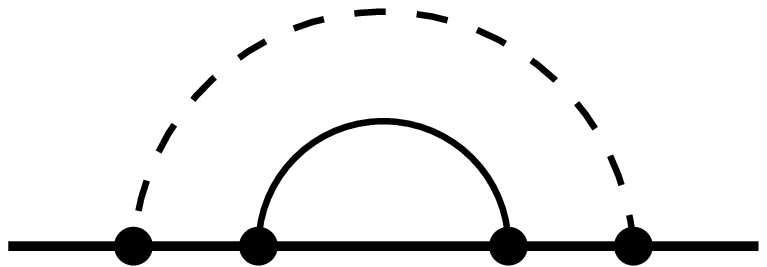}{0} & \jcd{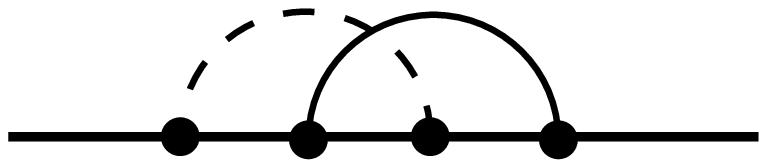}{2} \\
\jcd{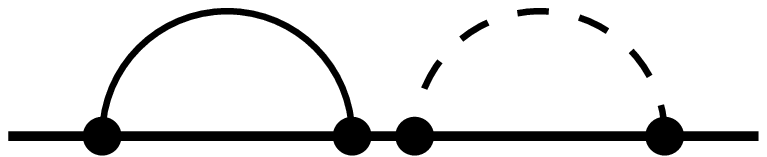}{0} & \jcd{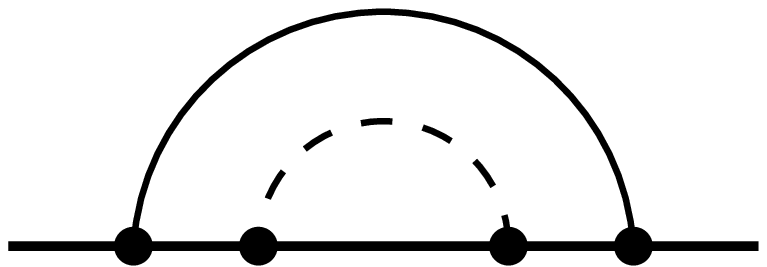}{4} & \jcd{pim}{2} \\
\jcd{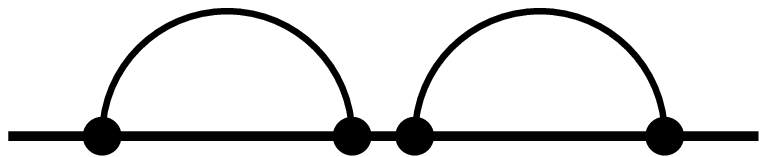}{0} & \jcd{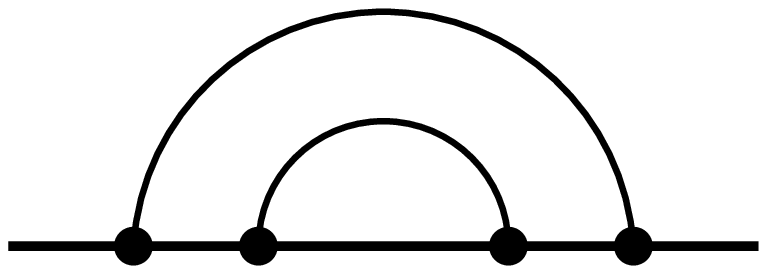}{0} & \jcd{pip}{2} \\
             &&\vspace{-5pt}\\
             \hline \\&&\vspace{-25pt}\\
\multicolumn{3}{|c|}{
 J^\pm_2(D^-_{m,n})=\left\{\begin{array}{r}
                        -6m-3,\mbox{\ for\ } m\geq0\\
                        2m-3,\mbox{\ for\ } m<0
                           \end{array}\right.
\qquad
 J^\pm_2(D^+_{m,n})=\left\{\begin{array}{r}
                        -6m-4,\mbox{\ for\ } m\geq0\\
                        2m-4,\mbox{\ for\ } m<0
                           \end{array}\right.}\\
             &&\vspace{-5pt}\\
             \hline \\&&\vspace{-25pt}\\
\multicolumn{3}{|c|}{\raisebox{8pt}{$J^\pm_2(D_n)=n$}
                     \makebox(0,22){}}\\
             \hline
\end{array}$$

\sbp{4.2.3.} Comparing this actuality table with the one of 3.4
one can notice that the $mod\ 2$ reduction of $J^\pm_2$ is an
invariant of the first order and $J^\pm_2 = J^-/2 (mod\ 2)$. This
proves the Corollary from {\bf 4.1}.

\sbp{4.2.4.} It would be interesting to find a Polyak--Viro style 
formula (see \cite{PV, CD}) for the invariant $J^\pm_2$ in terms of 
Gauss diagrams of the curve.

\sbp{4.2.5. Example.} Let us compute $J^\pm_2$ for the divide from 
{\bf 1.5.2} like it was done in {\bf 3.5}. The first 
step is pretty much the same
\def\jtpi#1{J^\pm_2\left(\ \mbox{\begin{picture}(30,10)(0,0)
        \put(0,-13){\epsfxsize=30pt \epsfbox{#1.eps}}
        \end{picture}}\ \right)}
$$\jtpi{ex34-0} = \jtpi{ex35-1} + \jtpi{ex35-2}$$
Now we compute separately the two terms of the right-hand side.
$$\begin{array}{ccl}
\jtpi{ex35-1} &=& \jtpi{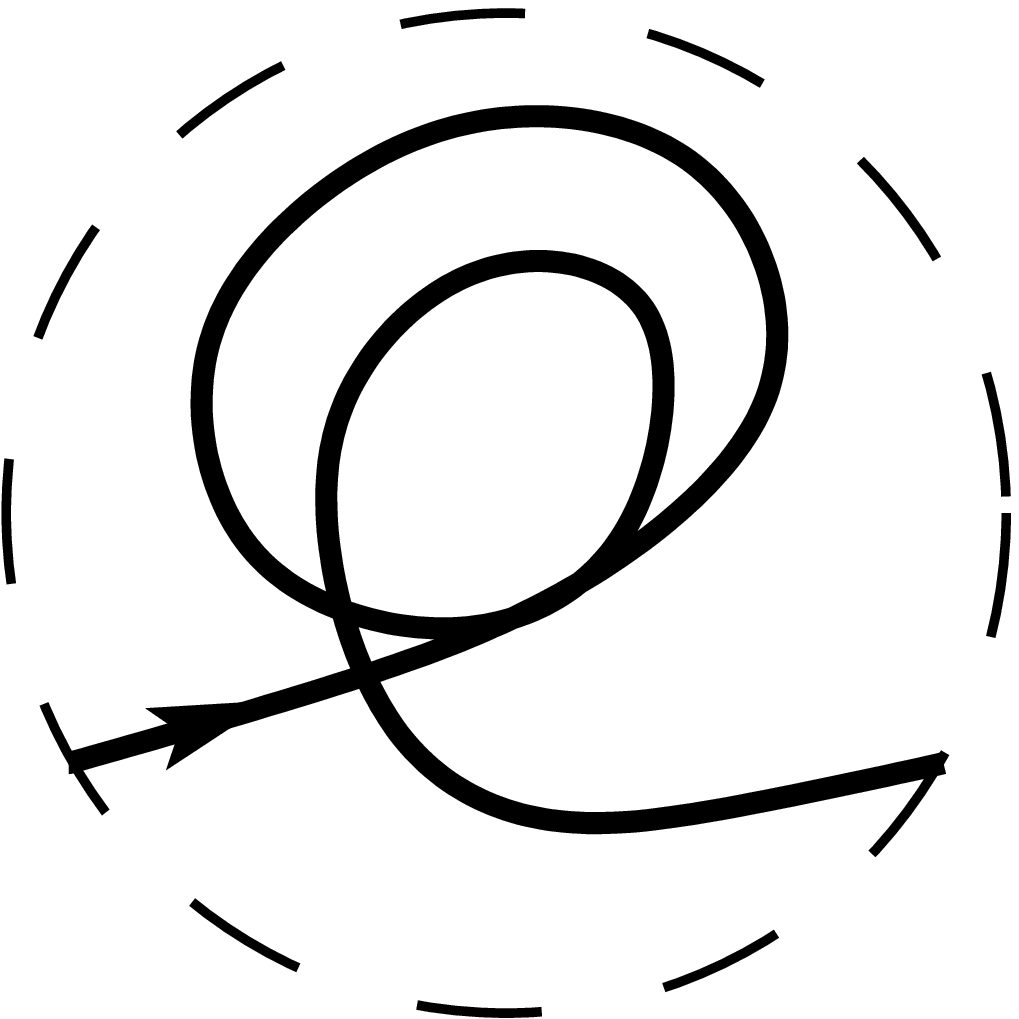} = \jtpi{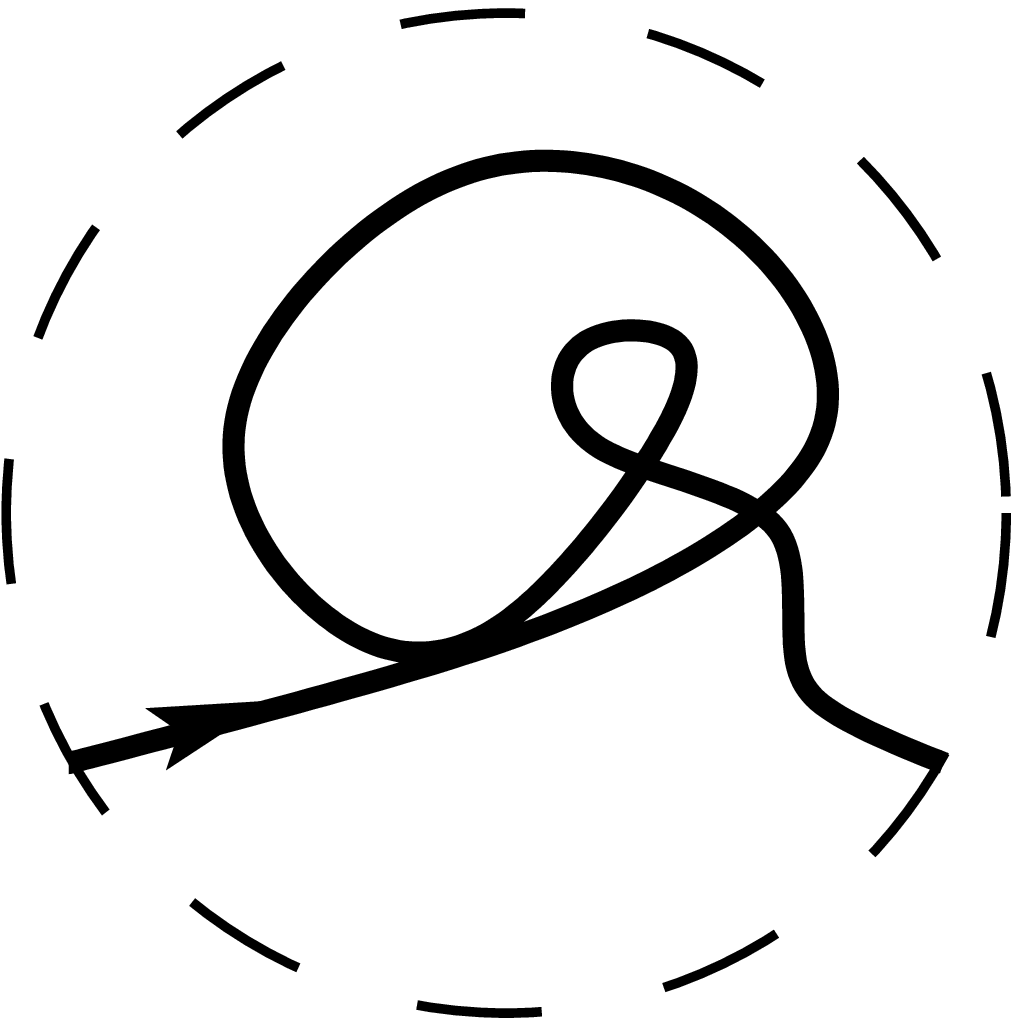} =
                  \jtpi{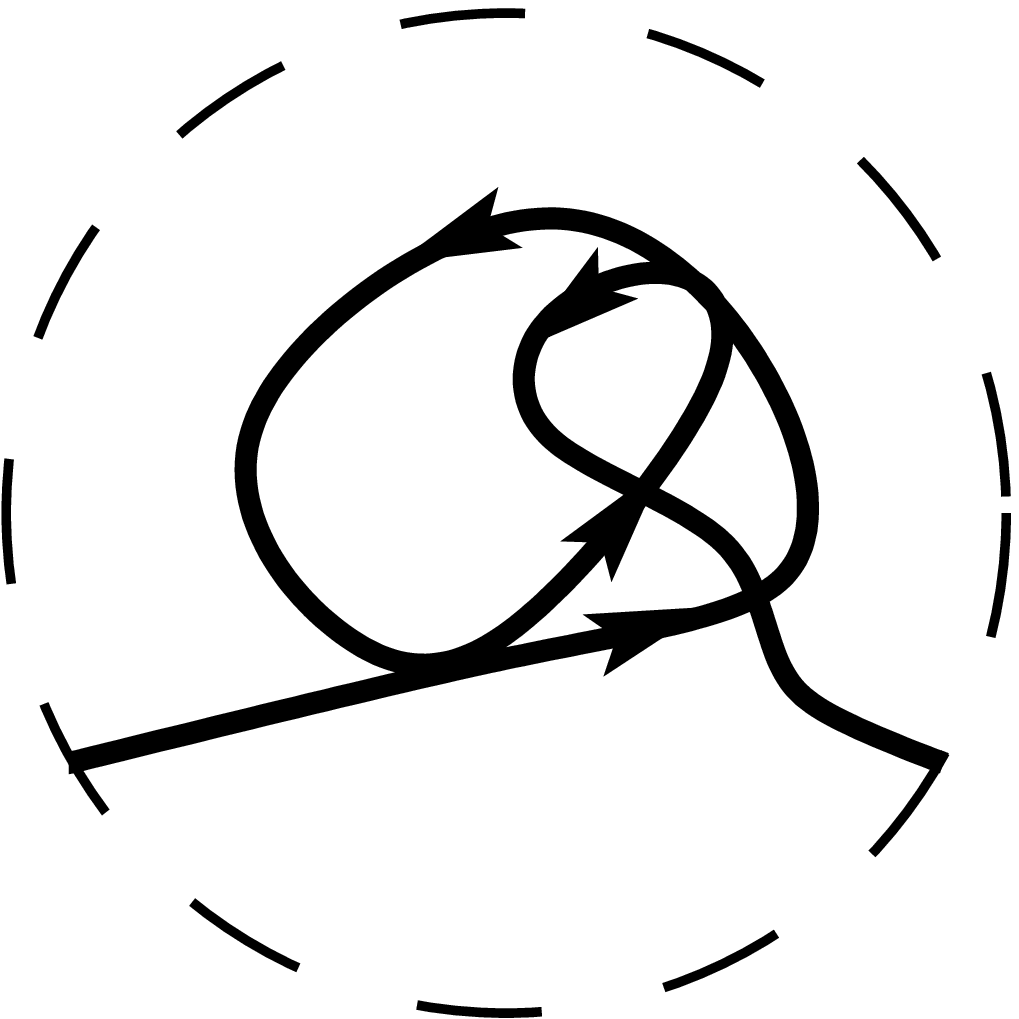} + \jtpi{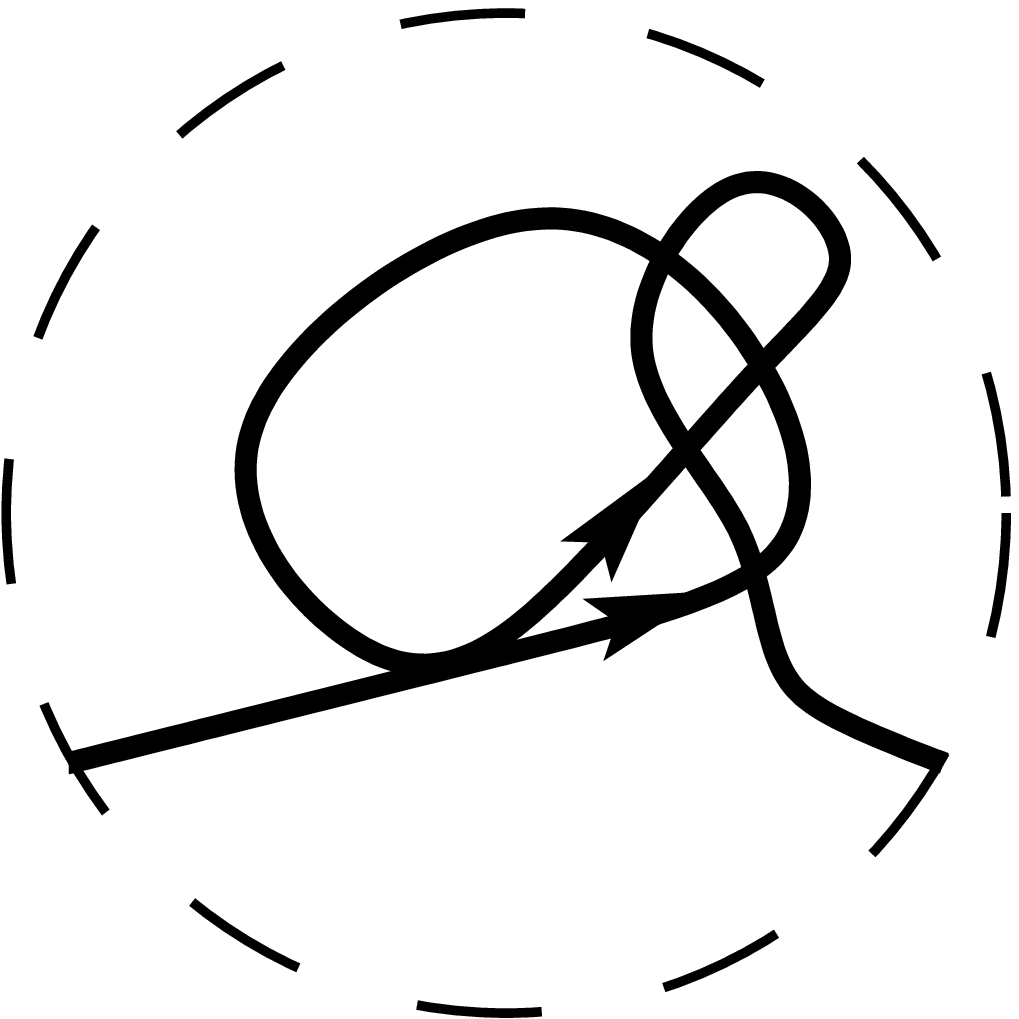}
                      \vspace{20pt}\\
&=& 
  J^\pm_2\Big(\ \begin{picture}(60,20)(0,0)
      \put(0,0){\epsfxsize=60pt \epsfbox{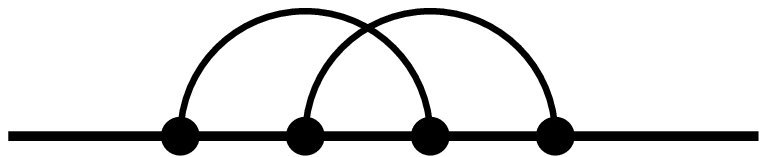}}
  \end{picture}\ \Big) + \jtpi{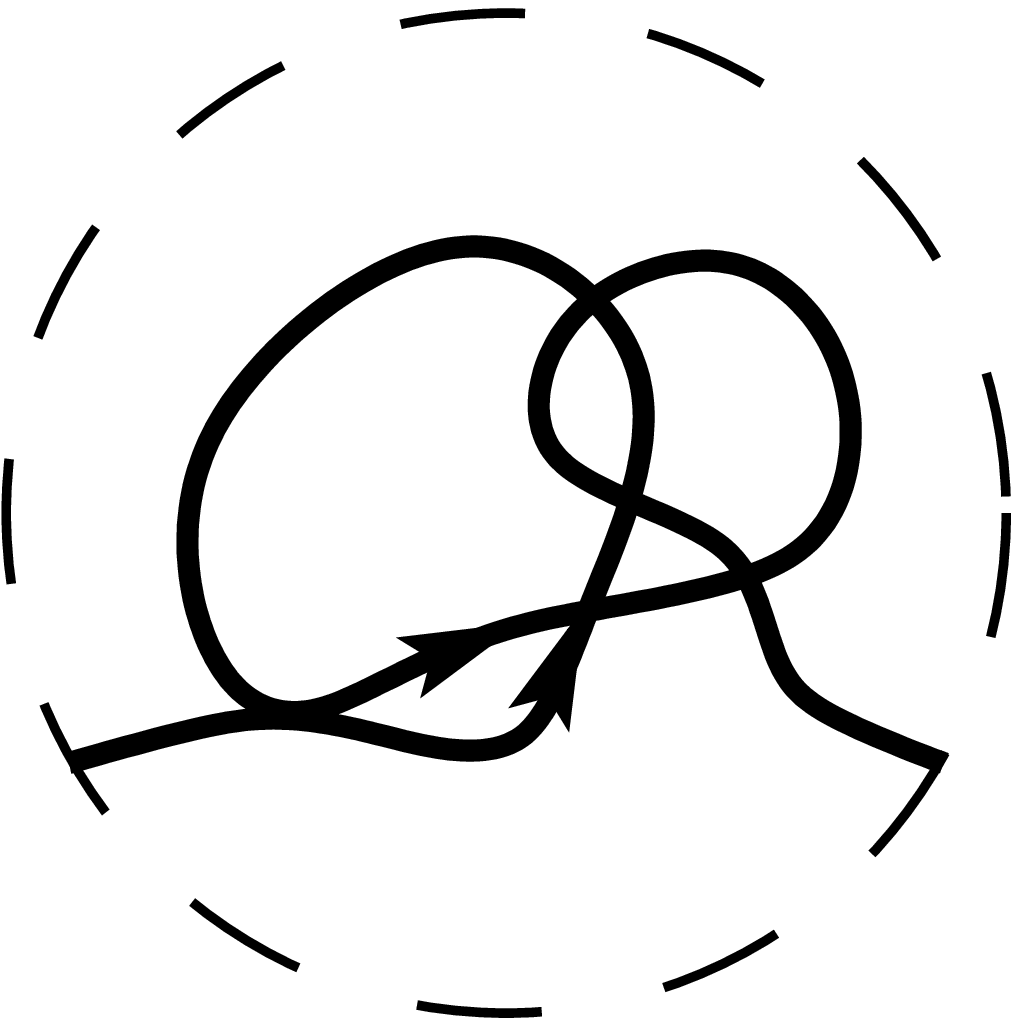} =
  2 - \jtpi{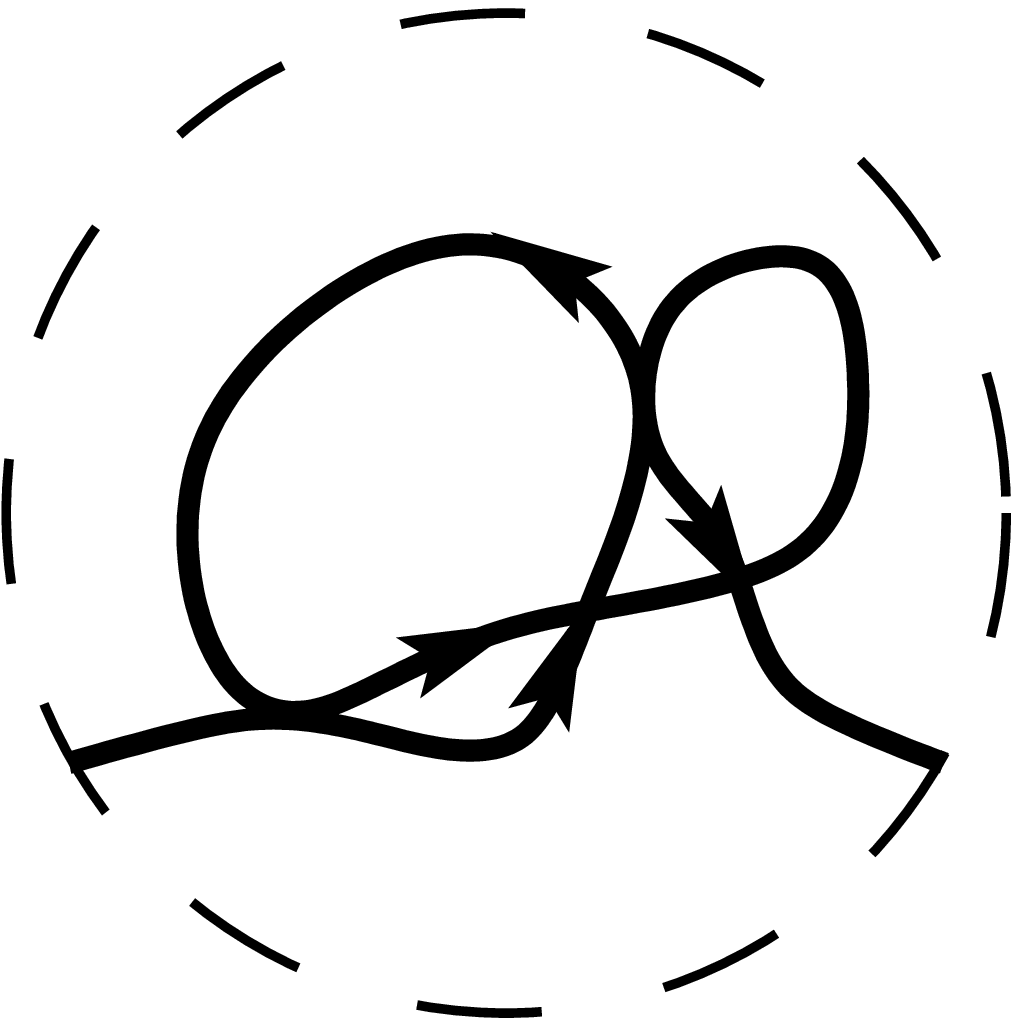} + \jtpi{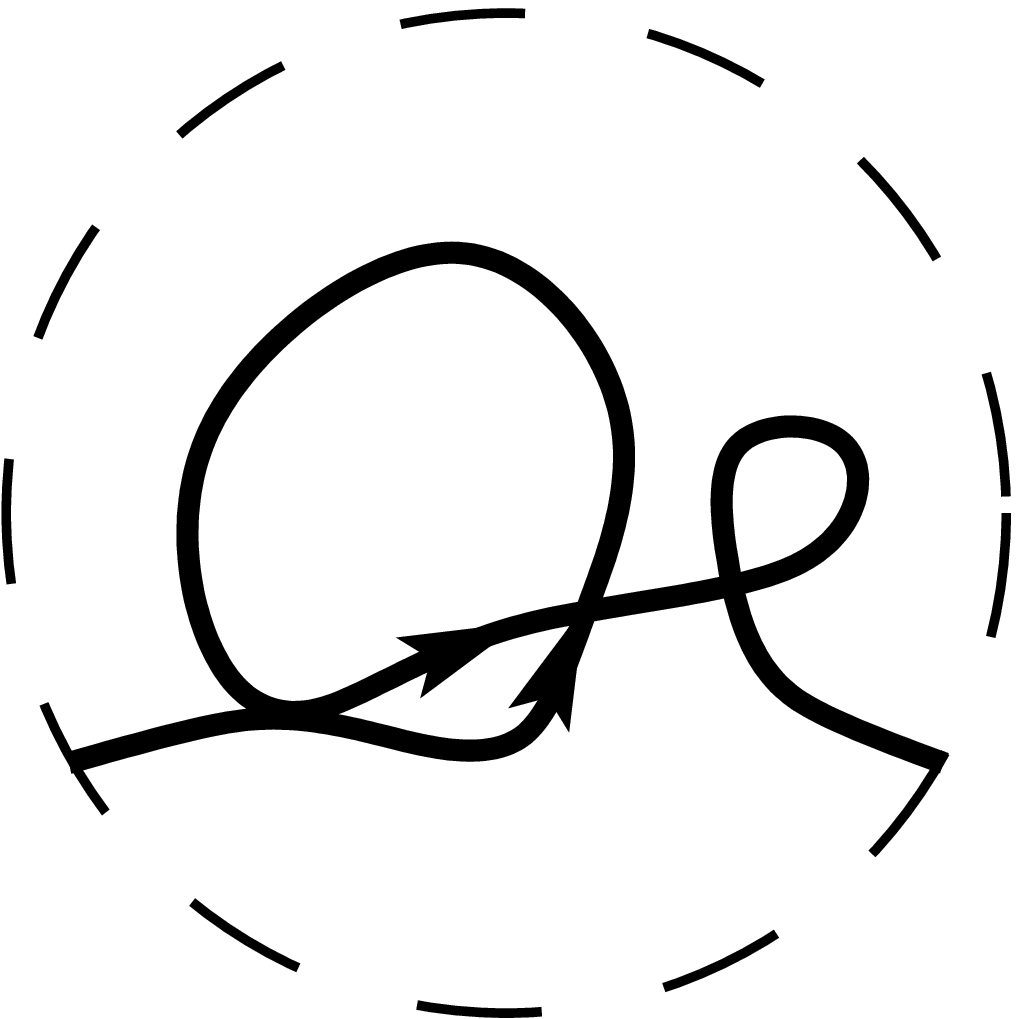}
                      \vspace{20pt}\\
&=&
  2-J^\pm_2\Big(\ \begin{picture}(60,20)(0,0)   
      \put(0,0){\epsfxsize=60pt \epsfbox{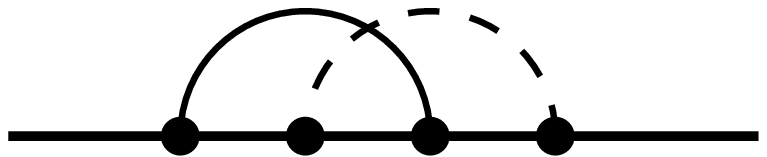}}
  \end{picture}\ \Big) + J^\pm_2(D^+_{0,1}) = -4
\end{array}
$$
For the second term we have
$$\begin{array}{ccl}
\jtpi{ex35-2} &=& \jtpi{ex35-3} = -\jtpi{ex35-4} + \jtpi{ex35-5}
                      \vspace{25pt}\\
&=&
   -\jtpi{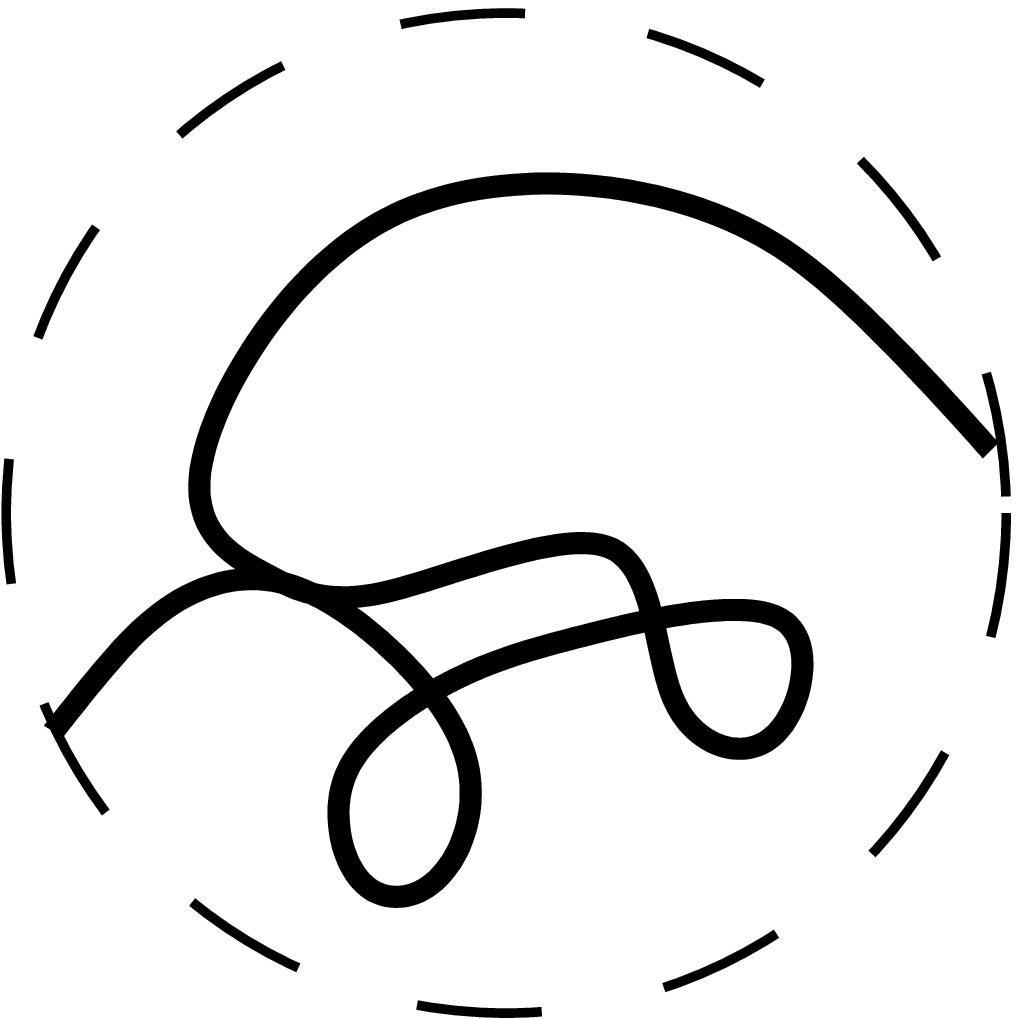} + J^\pm_2(D_2) = -J^\pm_2(D^-_{-2,0}) +2 =9
\end{array}\vspace{15pt}
$$
Totally, summing up, we get $\jtpi{ex34-0}=5$.

\sbp{4.2.6. Idea of a proof.} In the next section we describe how to
draw a diagram of $\L_D$ from the picture of the divide $D$. This
allows us to trace what happens with the knot $\L_D$ when the curve
$D$ changes by a direct (inverse) self-tangency move. Applying the
skein relation from {\bf 2.2} gives us the corresponding changings of
the Casson invariant. All this information can be summarized in the
actuality table of {\bf 4.2.2}.

\section{Hirasawa's Seifert surface for $\L_D$.}

\sbp{5.1.}
M.Hirasawa \cite{Hir} suggested a procedure to draw a picture of the
minimal genus Seifert surface for the link $\L_D$ and so to draw a
diagrams of the link $\L_D$. The other ways to draw diagrams see in
\cite{CP, Ch}. In this section we describe the Hirasawa's construction. 

First let us prepare the divide $D$ for drawing the Seifert surface as
follows. Choose an orientation of all branches of $D$. Then deform $D$
inside the unit disk so that:
\begin{itemize}
\item[$\bullet$] at every double point both branches are oriented
from left to right;
\item[$\bullet$] $x$ coordinates of double points and points of $D$
with vertical tangent are pairwise different.
\end{itemize}

Now to draw the Seifert surface we thicken every arc of $D$
to form a band and then modify the bands near double points and near
those points with vertical tangent where $D$ is oriented down as shown on
the pictures below.
\def\toto{\mbox{\Huge $\Longrightarrow$}}
$$\begin{picture}(250,80)(0,0)
  \put(0,15){\epsfxsize=50pt \epsfbox{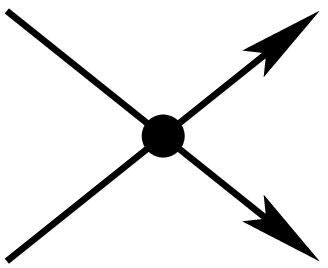}}
  \put(90,30){\toto}
  \put(140,0){\epsfxsize=100pt \epsfbox{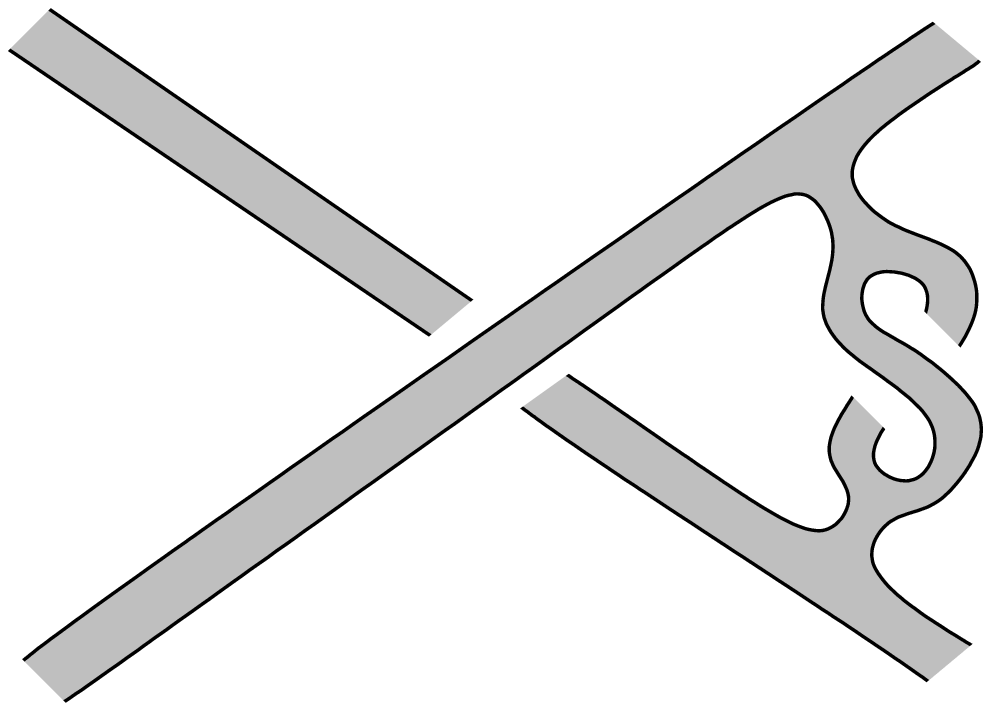}}
\end{picture}
$$

$$\begin{picture}(250,120)(0,0)
  \put(0,25){\epsfxsize=50pt \epsfbox{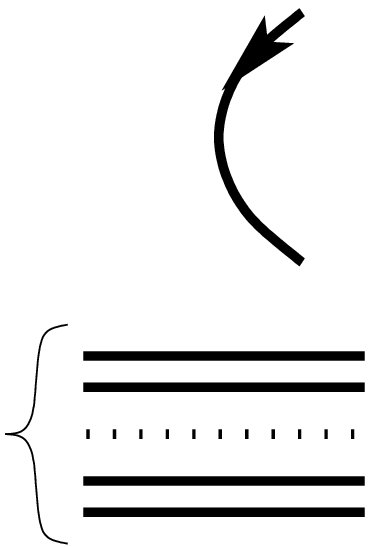}}
  \put(90,50){\toto}
  \put(140,0){\epsfxsize=100pt \epsfbox{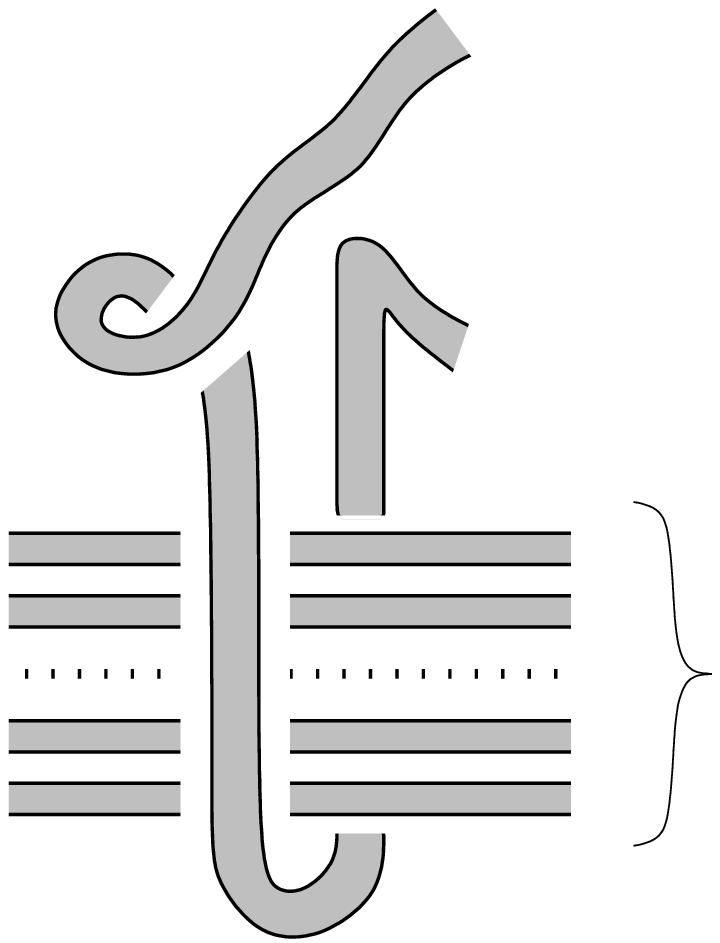}}
  \put(-100,40){\parbox{95pt}{\scriptsize
                   all arcs of $D$ located below the point with vertical
                   tangent line}}
  \put(245,35){\parbox{95pt}{\scriptsize
                   the bands corresponding to the  arcs of $D$
                   below the point with vertical
                   tangent line}}
\end{picture}
$$

$$\begin{picture}(250,120)(0,0)
  \put(0,25){\epsfxsize=50pt \epsfbox{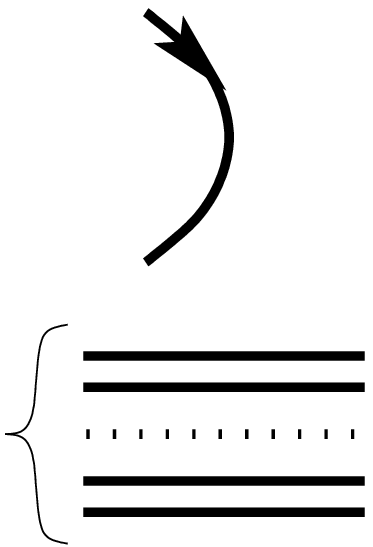}}
  \put(90,50){\toto}
  \put(140,0){\epsfxsize=100pt \epsfbox{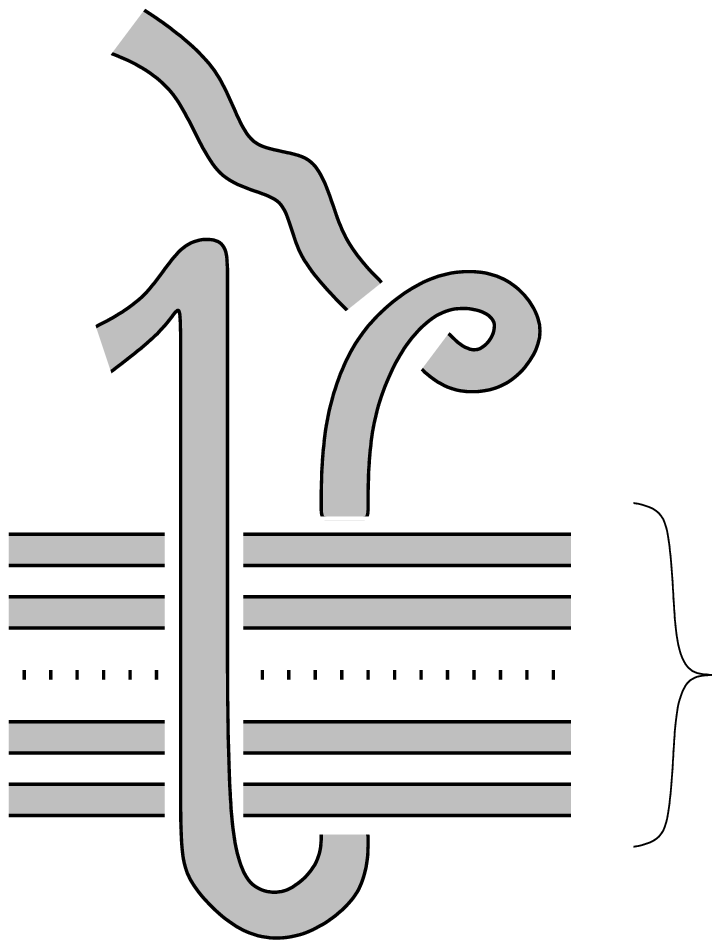}}
  \put(-100,40){\parbox{95pt}{\scriptsize
                   all arcs of $D$ located below the point with vertical
                   tangent line}}
  \put(245,35){\parbox{95pt}{\scriptsize
                   the bands corresponding to the  arcs of $D$
                   below the point with vertical
                   tangent line}}
\end{picture}
$$

\sbp{5.2. Example.} For the divide $D$ with a single curl (corresponding
to the singularity $A_2$) this construction gives the following Seifert
surface.
$$\begin{picture}(350,90)(0,0)
  \put(0,25){\epsfxsize=50pt \epsfbox{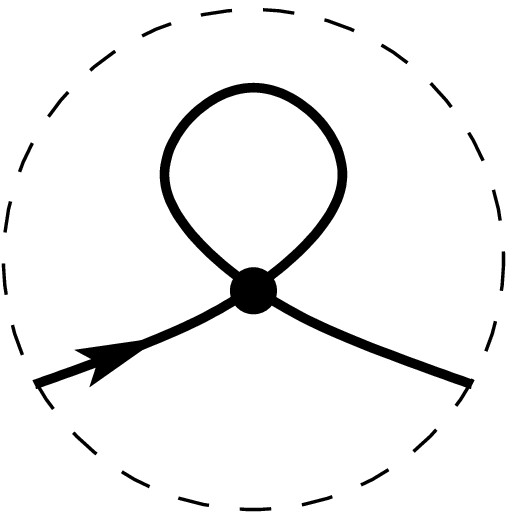}}
  \put(90,40){\toto}
  \put(140,0){\epsfxsize=100pt \epsfbox{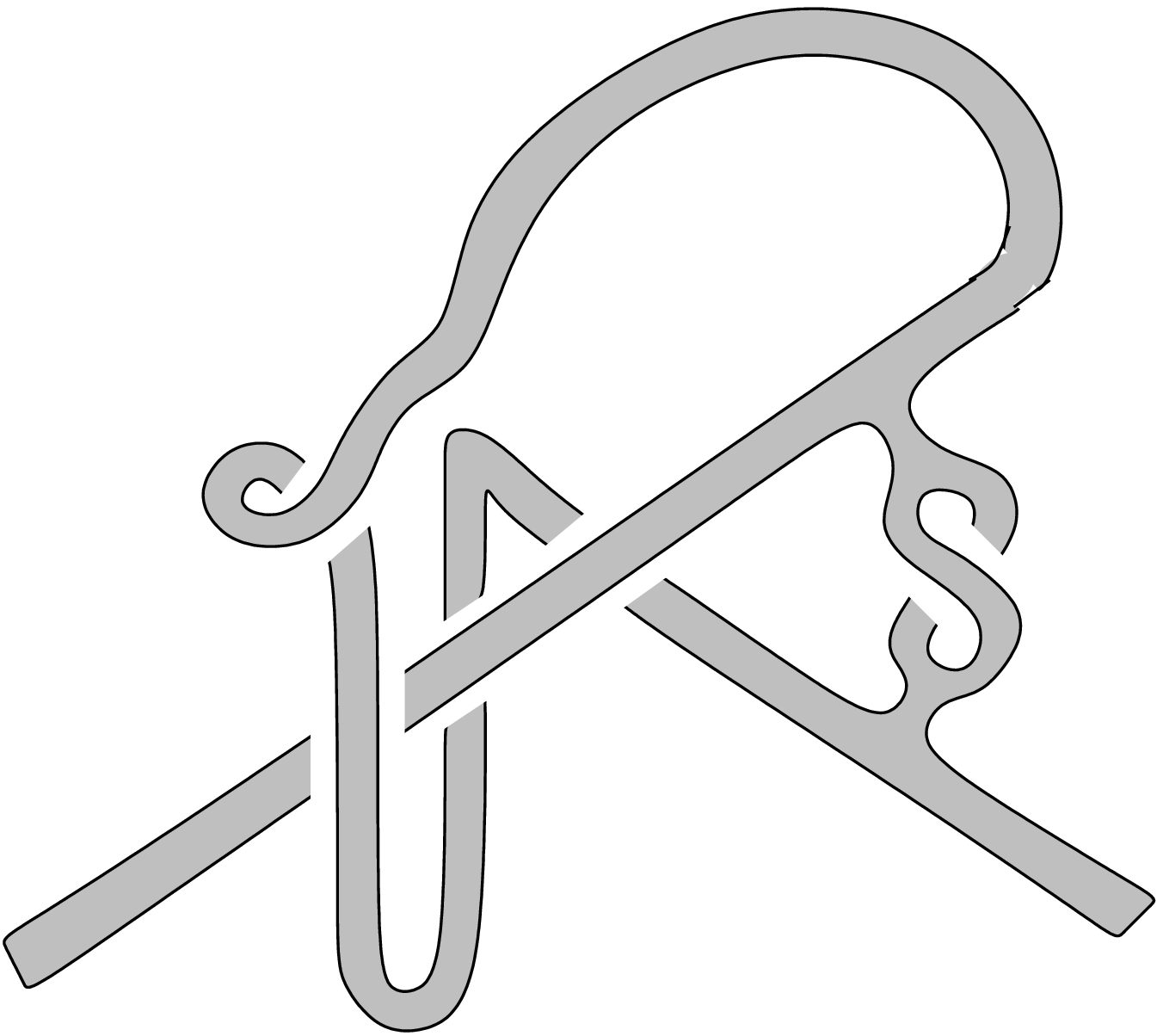}}
  \put(250,40){\mbox{\Large $\cong$}}
  \put(280,30){\epsfxsize=60pt \epsfbox{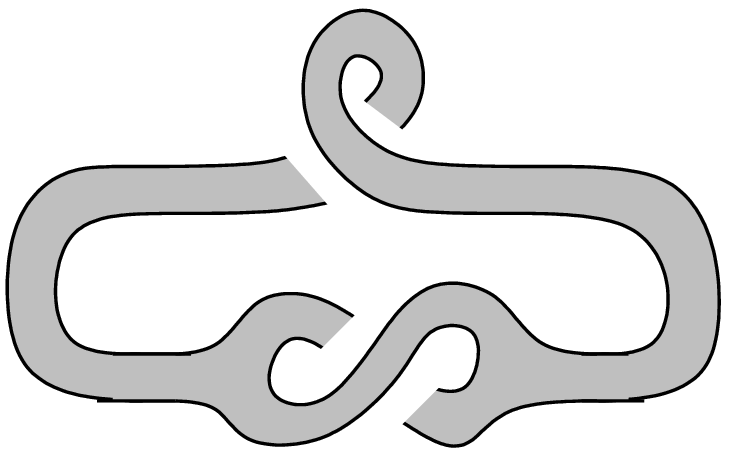}}
\end{picture}
$$
The corresponding knot $\L_D$ is the trefoil.
$$\begin{picture}(310,40)(0,0)
  \put(0,0){\epsfxsize=70pt \epsfbox{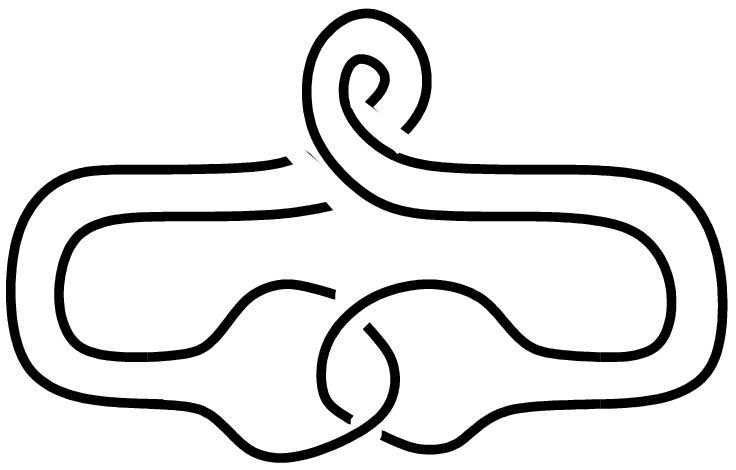}}
  \put(80,15){\mbox{\Large $\cong$}}
  \put(110,3){\epsfxsize=70pt \epsfbox{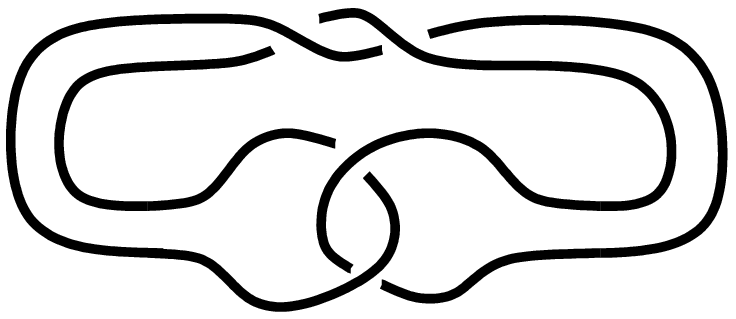}}
  \put(190,15){\mbox{\Large $\cong$}}
  \put(220,0){\epsfxsize=50pt \epsfbox{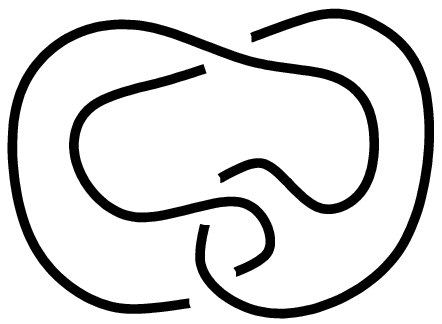}}
  \put(285,15){\mbox{\Large $=\ 3_1 $}}
\end{picture}
$$

\vskip 1truepc \noindent \it
 The Ohio State University,\\
Mansfield campus,\\
1680 University Drive,\\
Mansfield, OH 44906\\
{\rm E-mail:} {\tt chmutov@math.ohio-state.edu}


\begin{thebibliography}{BNG}

\bibitem[AC$_0$]{AC-0} N.~A'Campo,
    {\it Le groupe de monodromie du d\'eploiement des singularit\'es
    isol\'ees de courbes planes} I,
    Math. Ann., {\bf 213} (1975) 1--32.

\bibitem[AC$_1$]{AC-1} N.~A'Campo,
    {\it Real deformations and complex topology of plane curve
    singularities}, Ann. Fac. Sci. Toulouse Math. (6) 8 (1999),
    no. 1, 5--23 (see also {\tt alg-geom/9710023}).

\bibitem[AC$_2$]{AC-2} N.~A'Campo,
    {\it Generic immersions of curves, knots, monodromy and gordian
    number}, Inst. Hautes Etudes Sci. Publ. Math. No. 88 (1998),
    151--169 (1999) (see also {\tt math.GT/9803081}).

\bibitem[AC$_3$]{AC-3} N.~A'Campo,
    {\it Private communication}, January, 2001.

\bibitem[Ar]{Ar} V.~I.~Arnold, {\it Topological invariants of plane
        curves and caustics}, University Lecture Series, Vol. {\bf 5},
        AMS, Providence, RI, 1994.

\bibitem[AGV]{AGV} V.~I.~Arnold, S.~M.~Gusein-Zade, A.~N.~Varchenko,
    {\it Singularities of differentiable maps}, Vol. II,
    Birkh\"auser, Boston, MA, 1988.

\bibitem[Bu]{Bu} W.~Burau, {\it Kennzeichnung der Schlauchknoten},
    Abh. Math. Sem. Univ. Hamburg {\bf 9} (1932) 125--133.

\bibitem[CD]{CD} S.~Chmutov, S.~Duzhin, {\it  Explicit formulas
    for Arnold's generic curve invariants}, Arnold--Gelfand Mathematical
    Seminars, Birkh\"auser, Boston, MA,  (1997) 125-138.

\bibitem[Ch]{Ch} S.~Chmutov, {\it Diagrams of divide links}, To appear
     in Proc. AMS (see also\\ {\tt math.GT/0205329}).

\bibitem[CP]{CP} O.~Couture, B.~Perron,
    {\it Representative braids for links associated to plane immersed
    curves},
    Journal of Knot Theory and Its Ramifications, {\bf 9} (2000)
    1--30.

\bibitem[GZ]{GZ} S.~M.~Gusein-Zade,
    {\it Intersection matrices for certain singularities of functions
     of two variables},
    Functional Anal. and its Appl., {\bf 8} (1974) 10--13.

\bibitem[GZN]{GZN} S.~M.~Gusein-Zade, S.~M.~Natanzon,
    {\it The Arf-invariant and the Arnold invariants of plane curves},
    The Arnold-Gelfand mathematical seminars,
    Birkh\"auser, Boston, MA, (1997) 267--280 .

\bibitem[Hir]{Hir} M.~Hirasawa,
    {\it Visualization of A'Campo's fibered links and unknotting
     operations},
    to appear in Topology and its Applications.

\bibitem[Ka]{Ka} Louis H.~Kauffman, {\it On knots}, Annals of Math. Studies
    {\bf 115}, Princeton University Press (1987).

\bibitem[Mor]{Mor} H.Morton,
    {\it Private communication}, November, 2001.

\bibitem[PV]{PV} M.~Polyak, O.~Viro,
    {\it On the Casson Knot Invariant},  preprint {\tt math.GT/9903158},
    to appear in
    Journal of Knot Theory and Its Ramifications.

\bibitem[Ta]{Ta} S.~Tabachnikov,
    {\it Invariants of smooth triple point free plane curves},
    Journal of Knot Theory and Its Ramifications, {\bf 5} (1996)
    531--552.

\bibitem[Va]{Va} V. A. Vassiliev, {\it Cohomology of knot spaces},
     Theory of Singularities and Its Applications (ed. V. I. Arnold),   
     Advances in Soviet Math., {\bf 1} (1990) 23--69.

\end{thebibliography}
\end{document}